\documentclass[12pt]{amsart}
\usepackage{amsmath,amssymb,amsthm}
\usepackage{ifthen,verbatim,here}

\usepackage{mathrsfs}
\usepackage{yhmath}
\usepackage{booktabs}
\usepackage{underscore}
\usepackage[numbered,autolinebreaks,useliterate]{mcode}
\usepackage{longtable}
\usepackage{color}
\usepackage{url}
\usepackage{pdfpages}
\usepackage[english]{babel}
\usepackage{verbatim,here}
\usepackage[T1]{fontenc}
\usepackage{floatflt,graphicx,graphics}
\usepackage{a4wide}
\graphicspath{{figures}}

\numberwithin{equation}{section}
\nonstopmode
\theoremstyle{plain}

\newtheorem{theorem}[equation]{Theorem}

\newtheorem{lemma}[equation]{Lemma}

\theoremstyle{definition}

\newtheorem{nonsec}[equation]{}

\newtheorem{mysubsection}[equation]{}
{\qed\bigskip}

\newcounter{alphabet}

\newcommand{\be}{\begin{eqnarray}}
\newcommand{\ee}{\end{eqnarray}}
\newcommand{\ba}{\begin{array}}
\newcommand{\ea}{\end{array}}
\newcommand{\ben}{\begin{eqnarray*}}
\newcommand{\een}{\end{eqnarray*}}

\newcommand{\B}{\mathbb{B}}

\newcommand{\capa}{\mathrm{cap}\,}

\newcommand{{\tth}}{\mathrm{th}}
\newcommand{{\sh}}{\mathrm{sh}}
\newcommand{{\ch}}{\mathrm{ch}}

\newcommand{\E}{{\mathcal E}}

\newcommand{\K}{\mathcal{K}}

\renewcommand{\Im}{{\,\operatorname{Im}\,}}
\renewcommand{\Re}{{\,\operatorname{Re}\,}}

\newcommand {\M} {\mathsf{M}}

\renewcommand{\i}{\mathrm{i}}

\newcommand{\bI}{{\bf I}}
\newcommand{\bM}{{\bf M}}
\newcommand{\bN}{{\bf N}}
\renewcommand{\Im}{{ \rm Im}\,}
\renewcommand{\Re}{{ \rm Re}\,}
\newcommand{\CC}{{\mathbb C}}
\newcommand{\DD}{\mathbb{D}}

\font\fFt=eusm10 
\font\fFa=eusm7  
\font\fFp=eusm5  
\def\K{\mathchoice
{\hbox{\,\fFt K}}
{\hbox{\,\fFt K}}
{\hbox{\,\fFa K}}
{\hbox{\,\fFp K}}}
\font\fFt=eusm10 
\font\fFa=eusm7  
\font\fFp=eusm5  
\def\E{\mathchoice
{\hbox{\,\fFt E}}
{\hbox{\,\fFt E}}
{\hbox{\,\fFa E}}
{\hbox{\,\fFp E}}}

\newcounter{minutes}
\divide\time by 60
\newcounter{hours}
\multiply\time by 60
\addtocounter{minutes}{-\time}

\begin{document}

\bibliographystyle{amsplain}
\title[Circular arc polygons,  numerical conformal mappings]
{
Circular arc polygons,  numerical conformal mappings, and moduli of quadrilaterals 
}

\def\thefootnote{}
\footnotetext{
\texttt{\tiny File:~\jobname .tex,
          printed: \number\year-\number\month-\number\day,
          \thehours.\ifnum\theminutes<10{0}\fi\theminutes}
}
\makeatletter\def\thefootnote{\@arabic\c@footnote}\makeatother

\author[M. Nasser]{Mohamed Nasser}
\address{Department of Mathematics, Statistics and Physics, Qatar University, Doha, Qatar}
\email{mms.nasser@qu.edu.qa}
\author[O. Rainio]{Oona Rainio}
\address{Department of Mathematics and Statistics, University of Turku, FI-20014 Turku, Finland}
\email{ormrai@utu.fi}
\author[A. Rasila]{Antti Rasila}
\address{Guangdong Technion, Shantou, Guangdong 515063, China}
\email{antti.rasila@iki.fi; antti.rasila@gtiit.edu.cn}
\author[M. Vuorinen]{Matti Vuorinen}
\address{Department of Mathematics and Statistics, University of Turku, FI-20014 Turku, Finland}
\email{vuorinen@utu.fi}
\author[T. Wallace]{Terry Wallace}
\address{Guangdong Technion, Shantou, Guangdong 515063, China}
\email{wallaceterrywayne@gmail.com}
\author[H. Yu]{Hang Yu}
\address{School of Science, Zhejiang Sci-Tech University, Hangzhou 310018, China}
\email{yuhang_zstu@163.com}
\author[X. Zhang]{Xiaohui Zhang}
\address{School of Science, Zhejiang Sci-Tech University, Hangzhou 310018, China}
\email{xiaohui.zhang@zstu.edu.cn}

\begin{abstract}

We study numerical conformal mappings of planar Jordan domains with boundaries consisting of finitely many circular arcs and compute the moduli of quadrilaterals for these domains. Experimental error estimates are provided and, when possible, comparison to exact values or other methods are given. The main ingredients of the computation are boundary integral equations combined with the fast multipole method.   
\end{abstract}

\maketitle
\section{Introduction}
Many applications of complex analysis make use of conformal mappings or harmonic functions in an essential way. Examples include modeling of airfoils, analysis of turbulence, design of dams, and mathematical theory of electricity  \cite{af}. A first example is a mathematical model of a physical condenser, which is a planar domain $G$ together with a compact set $E\subset G$ \cite{ps}. The pair $(G,E)$ is called a \emph{condenser} and its \emph{capacity} is defined as \cite{du}
\begin{equation*} 
{\rm cap}(G,E)=\inf_{u\in A}\int_{G}|\nabla u|^2 dm,
\end{equation*}
where $A$ is the family of all harmonic functions with $u(x)\geq1$ for all $x\in E$ and $u(x)\to0$ when $x\to\partial G$. In concrete applications the sets $E$ and $\partial G$ have a simple geometry, both have a finite number of components, each component being a piecewise smooth curve. In this case, it is known that the infimum is attained by a harmonic function. A second example is  
to solve  the following Dirichlet-Neumann
boundary value problem for the Laplace equation in a planar Jordan domain $G$ with a piecewise smooth 
boundary $\partial
G = \cup_{k=1}^4\partial G_k$ where the sets  $\partial G_k$ occur in positive order---the sets $\partial G_k$ and
the domain $G$ define a \emph{quadrilateral}.
This problem is

\begin{equation*} 
\left\{\begin{matrix}
    \Delta u =&\ 0,\quad \text{on}\ {G,} \\
    u =&\ 1,\quad \text{on}\ {\partial G_1,}\\
    u =&\ 0,\quad \text{on}\ {\partial G_3,}\\
    \partial u/\partial n  =&\ 0,\quad \text{on}\ {\partial G_2,}\\
    \partial u/\partial n  =&\ 0,\quad \text{on}\ {\partial G_4.}\\
\end{matrix}\right.
\end{equation*}
If $u$ is a solution function to this problem, then the expression  $
\int\!\!\int_Q |\nabla u|^2 dxdy \,$ defines the \emph{ modulus}  of the quadrilateral.
Sometimes, conformal mappings can be applied  to simplify the geometry and, in fact, to solve
this Dirichlet-Neumann problem. By the Riemann mapping theorem
we know that a Jordan domain $G$ is conformally equivalent to the unit disk, but this theorem does
not give a clue for finding the conformal map. In the case of polygonal domains, widely used numerical methods
based on the Schwarz-Christoffel transformation  have been developed by D. Gaier \cite{Ga}, L.N. Trefethen and T.A. Driscoll \cite{dt}, and N. Papamichael and N. Stylianopoulos \cite{ps10}. 
Jordan domains where the boundary is a union of finitely many circular arcs, have 
been studied by
P. Brown and M. Porter \cite{bp}, \cite{p}, by U. Bauer and W. Lauf \cite{BL}, and, in particular, by D. Crowdy \cite{c}.  See also \cite{BjG,Ho,Tr}.
We give new numerical methods for this same case and present our results in the form of numerical tables, graphics, and analysis of algorithm performance. The method is based on boundary integral equations (R. Kress, \cite{kre14}) as developed and implemented by 
M.M.S. Nasser in a series of papers during the past two decades,
see e.g. \cite{KNH}, \cite{Nas-Siam1}-\cite{Nas-PlgCir}. The method uses the fast multipole method implementation from \cite{Gre-Gim12} for the speed-up of solving linear equations. In a recent series of papers (\cite{Nas-PlgCir}-\cite{nv1}), the method was applied for the capacity computation of planar condensers and for the study of isoperimetric problems for capacity. In particular, we will make use of the very recent results in \cite{Nas-PlgCir}.

The structure of this paper is as follows. In Section 2, we give the basic facts about the boundary integral method and its efficient implementation. Section 3 describes the method of Kress \cite{kre90} for the treatment of non-smooth boundary points in numerical integration and a refinement principle in trapezoid rule. Section 4 contains applications of the presented method to several circular arc polygonal quadrilaterals. Several numerical examples for the computation of the exterior modulus of quadrilaterals are presented in Section 5.
Finally, in Section 6, the proposed method is applied to gear domains.

\section{Preliminary notions}

The capacity of a condenser defined in the introduction can be
defined in many equivalent ways as shown in \cite{du}, \cite{GMP}. First, the family $A$ may be replaced by several other families by \cite[Lemma 5.21, p. 161]{GMP}. Furthermore,
\begin{align} \label{capmod}
\capa(G,E)=\M(\Delta(E,\partial G;G)),    
\end{align}
where $\Delta(E,\partial G;G)$ is the family of all curves joining $E$ with the boundary $\partial G$ in the domain $G$ and $\M$ stands for the modulus of a curve family \cite[Thm 5.23, p. 164]{GMP}. For the basic facts about capacities and moduli, the reader is referred to \cite{du, GMP, HKV}.

\begin{nonsec}{\bf Quadrilaterals.} \label{quad}
A Jordan domain in the complex plane $\mathbb C$ is a domain with boundary homeomorphic to the
unit circle. A {\it quadrilateral} is a Jordan domain $D$ together with four distinguished points $z_1,z_2,z_3,z_4 \in \partial D$ which define a positive orientation of the boundary. In other words, if we traverse the boundary,  then the points occur in the order of indices and the domain $D$ is on the left hand side. The quadrilateral
is denoted by $(D;z_1,z_2,z_3,z_4).$ The {\it modulus of the quadrilateral}
is a unique positive number $h$ such that $D$ can be conformally mapped by some conformal map $f$ onto the rectangle with vertices
$0,1,1+ih, ih$ such that 
$$f(z_1) = 0, \quad f(z_2) =1,  \quad f(z_3) =1+ih ,  \quad f(z_1) =ih\,.$$
The modulus is denoted ${\rm mod}(D;z_1,z_2,z_3,z_4).$  The following basic formula is often used:
\begin{equation}\label{reciprel}
{\rm mod}( D;z_1,z_2,z_3,z_4) \,{\rm mod}( D;z_2,z_3,z_4, z_1)=1\,.
\end{equation}
A simple example of a quadrilateral is the case when the domain $D$ is a rectangle with sides $a$ and $b$ and the points $z_1,z_2,z_3,z_4$ are its vertices. Depending on the labeling of the
vertices, the modulus is either $a/b$ or $b/a.$ An alternative equivalent definition is based on the Dirichlet-Neumann
problem mentioned in the introduction L.V. Ahlfors \cite[Thm~4.5, p.~63]{Ah}.
\end{nonsec}

\begin{nonsec}{\bf Quadrilateral modulus and curve families.} \label{quadrCurves}
The modulus of a quadrilateral $(D;z_1,z_2,z_3,z_4)$  is connected
with the modulus of the family of all curves in $D,$ joining the opposite boundary arcs $(z_2,z_3)$ and $(z_4,z_1),$ in a very simple way, as follows
\begin{equation} \label{2moduli}
{\rm mod}(D;z_1,z_2,z_3,z_4) = \M(\Delta((z_2,z_3), (z_4,z_1);D)) \,.
\end{equation}
\end{nonsec}

\begin{nonsec}{\bf Gr\"otzsch ring and elliptic integrals.} 
A {\it ring domain}, or briefly a {\it ring}, in the plane is a domain whose
complement has exactly two components. Such a domain can be conformally
mapped onto $\{z \in {\mathbb C}: r<|z|<1\}.$ The number $\log(1/r)$ is called {\it the modulus of the ring}. One can also consider a ring $R$ with
complementary components $E, F \subset {\mathbb C}$ as a condenser 
$({\mathbb C} \setminus E, F)$ and define the capacity of the ring as
\[
{\rm cap} R= {\rm cap} ({\mathbb C} \setminus E, F) =
2\pi/  {\rm mod}R .
\]
The formula \eqref{capmod} gives now the connection between the modulus
of a ring and the modulus of the family of curves joining its boundary components.
 The so called Gr\"otzsch ring $G_r={\mathbb B}^2 \setminus [0,r], 0<r<1,$ is frequent in the study of capacities.
 In the case $n=2, r\in(0,1),$ the following explicit
formulas hold for the capacity of this ring \cite[(7.18), p. 122]{HKV}, 
\begin{equation} \label{capGro}
{\rm cap} G_r=\frac{2\pi}{\mu(r)}\,; \quad \mu(r)=\frac{\pi}{2}\frac{\K'(r)}{\K(r)},
\end{equation}
where $\K(r)$ and $\K'(r)$ are the elliptic integrals of the first kind
\begin{equation} \label{eq:ellipK}
\K(r)=\int^1_0 \frac{dx}{\sqrt{(1-x^2)(1-r^2x^2)}},
\quad \K(r)=\K(r'), \quad r'=\sqrt{1-r^2}\,.
\end{equation}
The elliptic integrals $\E(r)$ and $\E'(r)$ of the second kind are 
\begin{equation} \label{eq:ellipE}
\E(r)=\int^{\pi/2}_0 \sqrt{1-r^2\sin^2 x}\,dx,
\quad \E'(r)=\E(r'), \quad r'=\sqrt{1-r^2}\,.
\end{equation}
\end{nonsec}

{\color{red} }
\begin{nonsec}{\bf History of numerical conformal mapping.}
The Schwarz-Christoffel method described in the introduction has a long history
which goes back to the nineteenth centure, see \cite{BG}. During the past fifty
years, the development of computational methods has revolutionized the applications
of numerical conformal mapping, see \cite{Ga, dt,ps10,ky}.
\end{nonsec}

\section{Boundary integral method}

\begin{mysubsection}{\bf The integral equation.}
We assume that the boundary $\Gamma = \partial G$ is a smooth Jordan curve parametrized by a $2\pi$-periodic  function $\eta : [0, 2\pi] \to \Gamma$ which is twice continuously differentiable and satisfies $\eta'(t) \neq 0$ for all $t\in[0,2\pi]$ (piecewise smooth boundaries will be considered in the next subsection).
The boundary $\Gamma$ is oriented such that $G$ is to the left of $\Gamma$, i.e., $\Gamma$ is oriented counterclockwise for bounded $G$ and clockwise for unbounded $G$. We denote by $H$ the space of all H\"older continuous real-valued functions on the boundary $\Gamma$.
\end{mysubsection}

Let $A : [0,2\pi] \to {\mathbb C}\setminus \{0\}$ be the complex function
\begin{equation}\label{eq:A}
A(t)= \begin{cases}
\eta(t)-\alpha, & \mbox{if $G$ is bounded}, \\
1, & \mbox{if $G$ is unbounded},
\end{cases}
\end{equation}
where $\alpha$ is a given auxiliary point in the domain $G$.
The generalized Neumann kernel $N(s,t)$ is defined by~\cite{WMN}
\begin{equation}\label{eq:N}
N(s,t) :=
\frac{1}{\pi}\Im\left(\frac{A(s)}{A(t)}\frac{\eta'(t)}{\eta(t)-\eta(s)}\right),\quad t \neq s,
\end{equation}
\begin{equation}\label{eq:Ntt}
N(t,t) :=
\frac{1}{\pi}\left(\frac{1}{2}\Im\frac{\eta''(t)}{\eta'(t)}-\Im\frac{A'(t)}{A(t)}\right).
\end{equation}
The kernel $N(s,t)$ is continuous on $[0,2\pi]\times [0,2\pi]$
 and hence, the integral operator $\bN$ defined on $H$ by
\[
\bN\gamma(s) := \int_{0}^{2\pi} N(s,t) \gamma(t) dt, \quad s\in [0,2\pi],
\]
is compact.  The integral equation involves also the kernel 
\begin{equation}\label{eq:M}
M(s,t) :=
\frac{1}{\pi}\Re\left(\frac{A(s)}{A(t)}\frac{\eta'(t)}{\eta(t)-\eta(s)}\right),  \quad s \neq t,
\end{equation}
which is singular and has the representation
\begin{equation}\label{eq:M1}
M(s,t) = -\frac{1}{2\pi}\cot\frac{s-t}{2}+M_1(s,t).
\end{equation}
Here the kernel $M_1$ is  continuous on $[0,2\pi]\times [0,2\pi]$ where
\begin{equation}\label{eq:M1tt}
M_1(t,t) :=
\frac{1}{\pi}\left(\frac{1}{2}\Re\frac{\eta''(t)}{\eta'(t)}-\Re\frac{A'(t)}{A(t)}\right).
\end{equation}
The integral operator $\bM$ defined on $H$ by
\[
\bM\gamma(s) := \int_{0}^{2\pi}  M(s,t) \gamma(t) dt, \quad s\in [0,2\pi],
\]
is singular, but is bounded on $H$~\cite{WMN}.

\begin{theorem}\label{thm:ie}
For a given function $\gamma\in H$, there exits a unique function $\rho\in H$ and a unique constant $h$ such that
the formula
\begin{equation}\label{eqn:Af}
\frac{\gamma+h+\i\rho}{A}
\end{equation}
defines the boundary values of an analytic function $f$ in $G$ with $f(\infty)=0$ for unbounded $G$. The function $\rho$ is the unique solution of the integral equation
\begin{equation}\label{eqn:ie}
(\bI - \bN) \rho = - \bM \gamma 
\end{equation}
and the constant $h$ is given by
\begin{equation} \label{eqn:h}
h = ( \bM \rho - (\bI - \bN) \gamma )/2.
\end{equation}
\end{theorem}

A MATLAB function \verb|fbie| for solving the integral equation~\eqref{eqn:ie} in $O(n\log n)$ operations, where $n$ is the number of nodes in the interval $[0,2\pi]$, is presented in~\cite{Nas-ETNA}. In the function \verb|fbie|, the integral equation is discretized using the Nystr\"om method with the trapezoidal rule and then using the GMRES method to solve the obtained linear system. The matrix-vector product in the GMRES method is computed using the MATLAB function $\mathtt{zfmm2dpart}$ in the toolbox $\mathtt{FMMLIB2D}$~\cite{Gre-Gim12}. Let \verb|et|, \verb|etp|, \verb|A|, \verb|gam|, and \verb|rho|, be $n\times1$ discretization vectors of the functions $\eta(t)$, $\eta'(t)$, $A(t)$, $\gamma(t)$, and $\rho(t)$, respectively. Then, the $n\times1$ vectors \verb|rho| and \verb|h| are computed by  
\[
[\verb|rho|,\verb|h|] = \verb|fbie|(\verb|et|,\verb|etp|,\verb|A|,\verb|gam|,\verb|n|,\verb|iprec|,\verb|restart|,\verb|gmrestol|,\verb|maxit|).
\]
Theoretically, all elements of the vector \verb|h| are equal to a constant $h$ in~\eqref{eqn:h}. Numerically, we will approximate the constant $h$ by the arithmetic mean of the elements of the vector \verb|h|. For the other parameters in $\verb|fbie|$, we choose $\mathtt{iprec}=5$, $\mathtt{gmrestol}=0.5\times 10^{-14}$, $\mathtt{restart}=[\,]$, and $\mathtt{maxit}=100$. This means that the tolerances of the methods are $0.5\times 10^{-15}$  for FMM and  $0.5\times 10^{-14}$ for GMRES. Moreover, the GMRES is used without restart, and the maximum number of GMRES iterations is $100$.  
The auxiliary points $\alpha$ in~\eqref{eq:A} need to be chosen sufficiently far away from the boundary $\Gamma$. For some domains (see Examples \ref{exm:Lp} and \ref{exp:Lcr} below), we need to choose $\alpha$ carefully to ensure the convergence of the method.

\begin{mysubsection}{ \bf Kress method.}

In this paper, we shall assume that $\Gamma$ is a piecewise smooth Jordan curve with a finite number of corner points such that each of these corner points is not a cusp.
We assume that the tangent vector of the boundary has only the first kind discontinuity at each corner point where the left tangent vector at each corner point is considered as the tangent vector at this point. 
For such boundaries $\Gamma$, the integral operator with the generalized Neumann kernel~\eqref{eqn:ie} is not compact, but this operator can be written as a sum of a compact operator and bounded non-compact operator with norm less than one in suitable function spaces~\cite{Nas-cr}. Hence, we can apply the Fredholm theory to the integral equation with the generalized Neumann kernel although the operator is not compact~\cite{kre90}. 

Using the method described above to solve the integral equation~\eqref{eqn:ie} when $\Gamma$ is a piecewise smooth Jordan curve yields only poor convergence since the solution of the integral equation~\eqref{eqn:ie} has a singularity in its first derivative in the vicinity of the corner points~\cite{kre90,Nas-cr}. 
To achieve a satisfactory accuracy, we first remove the discontinuity of the derivatives of the solution of the integral equation at the corner points by using a suitable substitution~\cite{kre90,kre91}. Then, the transformed equation can be solved using the above  method. 

Following Kress~\cite{kre90,kre91}, we define a bijective function $w : [0,2\pi] \to [0,2\pi]$ by
\begin{equation*}
w(t)=2\pi\frac{[v(t)]^p}{[v(t)]^p+[v(2\pi-t)]^p},
\end{equation*}
where 
\begin{equation*}
v(t)=\left(\frac{1}{p}-\frac{1}{2}\right)\left(\frac{\pi-t}{\pi}\right)^3+\frac{
1}{p}\frac{t-\pi}{\pi}+\frac{1}{2}.
\end{equation*}
The function $w$ is strictly monotonically increasing and infinitely differentiable function, and the integer $p\ge2$ is the grading parameter. In our numerical experiments below, we choose $p=3$.

Assume that the boundary $\Gamma$ has $m>0$ corner points and is parametrized by a $2\pi$-periodic function $\hat\eta(t)$ such that 
\begin{equation}\label{eq:corners}
\hat\eta(0),\quad \hat\eta(2\pi/m),\quad \hat\eta(4\pi/m), \quad\ldots, \quad \hat\eta(2(m-1)\pi/m)
\end{equation}
are the corner points of $\Gamma$. Assume also that $\hat\eta(t)$ is twice continuously differentiable with $\hat\eta'(t) \neq 0$ for all $t\in[0,2\pi]\backslash\{0,2\pi/m,\ldots,2(m-1)\pi/m\}$. 
Then we define a bijective, strictly monotonically increasing and infinitely 
differentiable function, $\delta:[0,2\pi]\to [0,2\pi]$, by~\cite{LSN17}
\[
\delta(t)= \begin{cases} 
 w(m t)/m,                         & t \in [0,2\pi/m), \\ 
[w(m t-2\pi)+2\pi]/m,              & t \in [2\pi/m,4\pi/m), \\ 
\quad \vdots & \\
[w(m t-2(m-1)\pi)+2(m-1)\pi]/m,    & t \in [2(m-1)\pi/m,2\pi].  
\end{cases}
\]
The function $\delta$ is at least $p$ times continuously differentiable since the function $w$ has a zero of order $p$ at the endpoints $t=0$ and $t=2\pi$~\cite[Theorem~2.1]{kre90}.

Then, we parametrize the boundary $\Gamma$ by $\eta(t)=\hat\eta(\delta(t))$ and hence $\eta'(t)=\hat\eta'(\delta(t))\delta'(t)$, $t\in[0,2\pi]$. With the new parametrization, the integral equation~\eqref{eqn:ie} is solved accurately using the above method as in the case of domains with smooth boundaries.

\end{mysubsection}

\begin{mysubsection}{\bf Parametrizing the boundaries of polygonal domains.}
In this paper, we assume that the boundary $\Gamma$ is a polygon consisting of a finite number of finite segments, circular arcs, or both  segments and circular arcs. The boundary $\Gamma$ will be parametrized as described in the preceding subsection. 
We discretize the interval $[0,2\pi]$ by $n$ equidistant nodes, 
\begin{equation}\label{eq:t_k}
t_k = (k-1) \frac{2 \pi}{n}, \quad k = 1, \ldots, n.
\end{equation}
Then the parametrization of the boundary $\Gamma$ is discretized by $\eta(s_1),\eta(s_s),\ldots,\eta(s_n)$ where $s_k=\delta(t_k)$, $k=1,2,\ldots,n$. Similarly, the derivative of the parametrization of the boundary is discretized by $\eta'(s_k)s'_k$ where $s'_k=\delta'(t_k)$, $k=1,2,\ldots,n$.  
A MATLAB function \verb|plgsegcirarcp.m| for computing such a parametrization can be downloaded from \url{https://github.com/mmsnasser/circa}.
To use this MATLAB function, assume that $\Gamma$ is a polygon with $m$ finite vertices $v_1,v_2,\ldots,v_m$ with $v_{m+1}=v_1$. 
For the parametrization $\eta(t)$ of the boundary $\Gamma$, it follows from~\eqref{eq:corners} and~\eqref{eq:t_k} that
\begin{equation}\label{eq:v_k}
v_k=\eta\left((k-1)\frac{2\pi}{m}\right)=\eta\left(t_{(k-1)n/m+1}\right), 
\quad k=1,2,\ldots,m.
\end{equation}
For $k=1,2,\ldots,m$, we assume that the center of the portion of the boundary between $v_k$ and $v_{k+1}$ is $c_k$ if the portion is a circular arc and $c_k=\infty$ if the portion is a segment. 
Further, if the portion of the boundary between $v_k$ and $v_{k+1}$ is a circular arc, then we introduce an indicator $d_k$ for it with value $1$ if the arc is positively oriented with respect to the pertaining circle center $c_k$ and $-1$ otherwise. For segment portions of the boundary this indicator is $d_k=0$.
We assume here that $n$ is an integer multiple of $m$ so that each side of the polygon will be discretized by $n/m$ points. 
Define the vectors
\[
{\tt v=[v_1,\ldots,v_n], \quad c=[c_1,\ldots,c_n], \quad d=[d_1,\ldots,d_n],  
}
\]
Then, discretizations of the parametrization of the boundary and its derivative can be computed using the MATLAB function \verb|plgsegcirarcp.m| by calling
\[
{\tt [et,etp]=plgsegcirarcp(v,c,d,n/m)}.
\]

\end{mysubsection}

\begin{mysubsection}{\bf Conformal mapping onto the unit disk.}\label{sec:map}
In this subsection, we review a numerical method for the computation of the conformal mapping $w=\Phi(z)$ from a polygonal domain $G$ onto the unit disk $\DD=\Phi(G)$~\cite{nvs,Nas-cmft15}. 

For a bounded domain $G$, let $\gamma=-\log|\eta(t)-\alpha|$,  let $\rho$ be the unique solution of the integral equation~\eqref{eqn:ie}, and let the constant $h$ be given by~\eqref{eqn:h}. 
Then, the mapping function $\Phi$ with normalization
\begin{equation}\label{eq:Phi-cond}
\Phi(\alpha)=0, \quad \Phi'(\alpha)>0
\end{equation}
can be written for $z\in G\cup\Gamma$ as
\begin{equation}\label{eq:Phi-b}
\Phi(z)=c(z-\alpha)e^{(z-\alpha)f(z)}
\end{equation}
where $c=\Phi'(\alpha)>0$ and the function $f(z)$ is analytic in $G$ with the boundary values $A(t)f(\eta(t))=\gamma(t)+h+\i\rho(t)$. The boundary values of the conformal mapping are  now  given by
\[
\Phi(\eta(t))=c(\eta(t)-\alpha)e^{\gamma(t)+h+\i\rho(t)}.
\]
Since $|\Phi(\eta(t))|=1$, it follows that $\Phi(\eta(t))$ can be written as
\begin{equation}\label{eq:S}
\Phi(\eta(t))=e^{\i S(t)}, \quad t\in[0,2\pi],
\end{equation}
where 
\begin{equation}\label{eq:S-b}
S(t)=\arg(\eta(t)-\alpha)+\rho(t)
\end{equation}
is a one-to-one increasing function on $[0,2\pi]$ with $S(2\pi)-S(0)=2\pi$. The function $S(t)$ is known as the boundary corresponding function~\cite[p.~380]{hen}. Differentiation yields
\[
S'(t)=\Im\left[\frac{\eta'(t)}{\eta(t)-\alpha}\right]+\rho'(t).
\]

When $G$ is an unbounded domain, we assume $\gamma=\log|\eta(t)-z_1|$ where $z_1$ is a given point in the exterior domain of $G$, i.e., $z_1$ is in the bounded domain enclosed by $\Gamma$. We assume also that $\rho$ is the unique solution of the integral equation~\eqref{eqn:ie} and the constant $h$ is given by~\eqref{eqn:h}. 
Then, the mapping function $\Phi$ with the normalization
\[
\Phi(\infty)=0, \quad \lim_{z\to\infty}\left(z\Phi(z)\right)>0
\]
can be written for $z\in G\cup\Gamma$ as
\begin{equation}\label{eq:Phi-u}
\Phi(z)=\frac{c}{z-z_1}e^{f(z)}
\end{equation}
where $c=\lim_{z\to\infty}\left(z\Phi(z)\right)>0$ and the function $f(z)$ is analytic in $G$ with the boundary values $A(t)f(\eta(t))=\gamma(t)+h+\i\rho(t)$ and $f(\infty)=0$. Note that $A(t)=1$ for unbounded domains. Hence, the boundary values of the conformal mapping are given by
\[
\Phi(\eta(t))=\frac{c}{\eta(t)-z_1}e^{\gamma(t)+h+\i\rho(t)}.
\]
The mapping function $\Phi$ maps the boundary $\Gamma$ onto the unit circle $\partial\DD=\Phi(\Gamma)$. Although the orientation of $\Gamma$ is clockwise ($G$ is on the left of $\Gamma$), the orientation of $\partial\DD$ will be counterclockwise. Then $|\Phi(\eta(t))|=1$ and hence $\Phi(\eta(t))$ can be written as in~\eqref{eq:S}
where 
\begin{equation}\label{eq:S-u}
S(t)=-\arg(\eta(t)-z_1)+\rho(t)
\end{equation}
is a one-to-one increasing function on $[0,2\pi]$ with $S(2\pi)-S(0)=2\pi$. By differentiation, we obtain
\[
S'(t)=-\Im\left[\frac{\eta'(t)}{\eta(t)-z_1}\right]+\rho'(t).
\]

For both bounded and unbounded $G$, the function $\rho(t)$ is computed by solving the integral equation~\eqref{eqn:ie} and its derivative $\rho'(t)$ is computed by approximating $\rho(t)$ by a trigonometric interpolating polynomial and then differentiating the interpolating polynomial. This polynomial can be computed using FFT~\cite{Weg05}. Then the boundary corresponding function $S(t)$ is computed through~\eqref{eq:S-b} for bounded $G$ and by~\eqref{eq:S-u} for unbounded $G$. The boundary values of the mapping function $\Phi$ can be computed through~\eqref{eq:S}. Thus, the function
\begin{equation}\label{eq:zet}
\zeta(t)=\Phi(\eta(t))=e^{\i S(t)}, \quad t\in[0,2\pi], 
\end{equation}
is a parametrization of the unit circle. The values of the mapping function $w=\Phi(z)$ for $z\in D$ as well as the values of the inverse mapping function $z=\Phi^{-1}(w)$ for $w\in \DD$ can be computed using the Cauchy integral formula. For the direct mapping, we have
\[
\Phi(z)=\frac{1}{2\pi\i}\int_{\Gamma} \frac{\Phi(\eta)}{\eta-z}d\eta
=\frac{1}{2\pi\i}\int_{0}^{2\pi} \frac{\Phi(\eta(t))}{\eta(t)-z}\eta'(t)dt 
=\frac{1}{2\pi\i}\int_{0}^{2\pi} \frac{\zeta(t)}{\eta(t)-z}\eta'(t)dt, 
\quad z\in D.
\]
For the inverse mapping, we have 
\[
\Phi^{-1}(w)=\frac{1}{2\pi\i}\int_{\partial\DD} \frac{\Phi^{-1}(\zeta)}{\zeta-w}d\zeta
=\frac{1}{2\pi\i}\int_{0}^{2\pi} \frac{\Phi^{-1}(\zeta(t))}{\zeta(t)-w}\zeta'(t)dt 
=\frac{1}{2\pi\i}\int_{0}^{2\pi} \frac{\eta(t)}{\zeta(t)-w}\zeta'(t)dt, 
\quad w\in \DD.
\]

The mapping function $\Phi$ maps the vertices $v_k$, $k=1,2,\ldots,m$, of the circular arc polygon domain $D$ onto points $w_k=\Phi(v_k)$, $k=1,2,\ldots,m$, on the unit circle. In literature, the points $w_k$, $k=1,2,\ldots,m$ are know as the preimages of the vertices $v_k$, $k=1,2,\ldots,m$. In view of~\eqref{eq:v_k}, it follows from~\eqref{eq:zet} that the preimages are given by
\[
w_k=\zeta(t_{(k-1)n/m+1})=e^{\i S(t_{(k-1)n/m+1})}, \quad k=1,2,\ldots,m.
\]

The implementation of the method presented in this subsection is given in the following MATLAB function \verb|mapdisk.m|.

\begin{lstlisting}
function [zet,zetp,c,S,Sp] = mapdisk(et,etp,n,zz,type)
if type=='b' 
    A = et-zz; k = 1;
elseif type=='u' 
    A = ones(size(et)); k = -1;
end
gam     = -k*log(abs(et-zz));
[rho,h] =  fbie(et,etp,A,gam,n,5,[],1e-14,200);
c       =  exp(-mean(h));
S       =  k*carg(et-zz)+rho;
Sp      =  k*imag(etp./(et-zz))+derfft(rho);
zet     =  exp(i*S);  zetp  =  exp(i*S).*Sp;
end
\end{lstlisting}

All computer codes of the computations presented in this paper can be found in the link \url{https://github.com/mmsnasser/circa}.

\end{mysubsection}

\begin{mysubsection}{\bf Modulus of quadrilaterals.}
Consider the quadrilateral $(D;z_1,z_2,z_3,z_4)$ where $D$ is a bounded simply connected polygonal domain and $z_1,z_2,z_3,z_4$ are four distinguished points on $\partial D$ with counterclockwise orientation. 
The modulus ${\rm mod}(D;z_1,z_2,z_3,z_4)$ can be computed in two steps. 
In the first step, we map the domain $D$ using the conformal mapping $w=\Phi(z)$ described in \S~\ref{sec:map} onto the unit disk $\DD$. 
The boundary $\partial D$ is then mapped onto the unit circle. 
Here, the points $z_1,z_2,z_3,z_4$ need not to be vertices of the polygon.
Assume that $z_k=\eta(\hat t_k)$ where $\hat t_k\in[0,2\pi]$, $k=1,2,3,4$. Then the four points $z_1,z_2,z_3,z_4$ will be mapped onto four points $w_1,w_2,w_3,w_4$ on the unit circle where, by~\eqref{eq:S}, $w_k=e^{\i S(\hat t_k)}$, $k=1,2,3,4$. By the conformal invariance of the modulus, we have
\[
{\rm mod}(D;z_1,z_2,z_3,z_4) = {\rm mod}(\DD;w_1,w_2,w_3,w_4).
\]

In the second step, the modulus ${\rm mod}(\DD;w_1,w_2,w_3,w_4)$ will be computed using the exact formula~\cite[(2.6.1)]{ps10}
\begin{equation}\label{eq:mod-disk}
{\rm mod}(\DD;w_1,w_2,w_3,w_4) = \frac{2}{\pi}\mu\left(1/\sqrt{k}\right), \quad k=|w_1,w_2,w_3,w_4|,
\end{equation}
where the absolute (cross) ratio $|w_1,w_2,w_3,w_4|$ is defined by~\cite[p.~33]{HKV}
\begin{equation}\label{eq:ration}
|w_1,w_2,w_3,w_4| = \frac{|w_1-w_3||w_2-w_4|}{|w_1-w_2||w_2-w_4|}.
\end{equation}
(Note that the definition of the cross-ratio here is different from the definition in~\cite[(1.10.5)]{ps10}.)

\begin{figure}[H]
\centerline{
\scalebox{0.55}{\includegraphics[trim=0cm 0cm 0cm 0cm,clip]{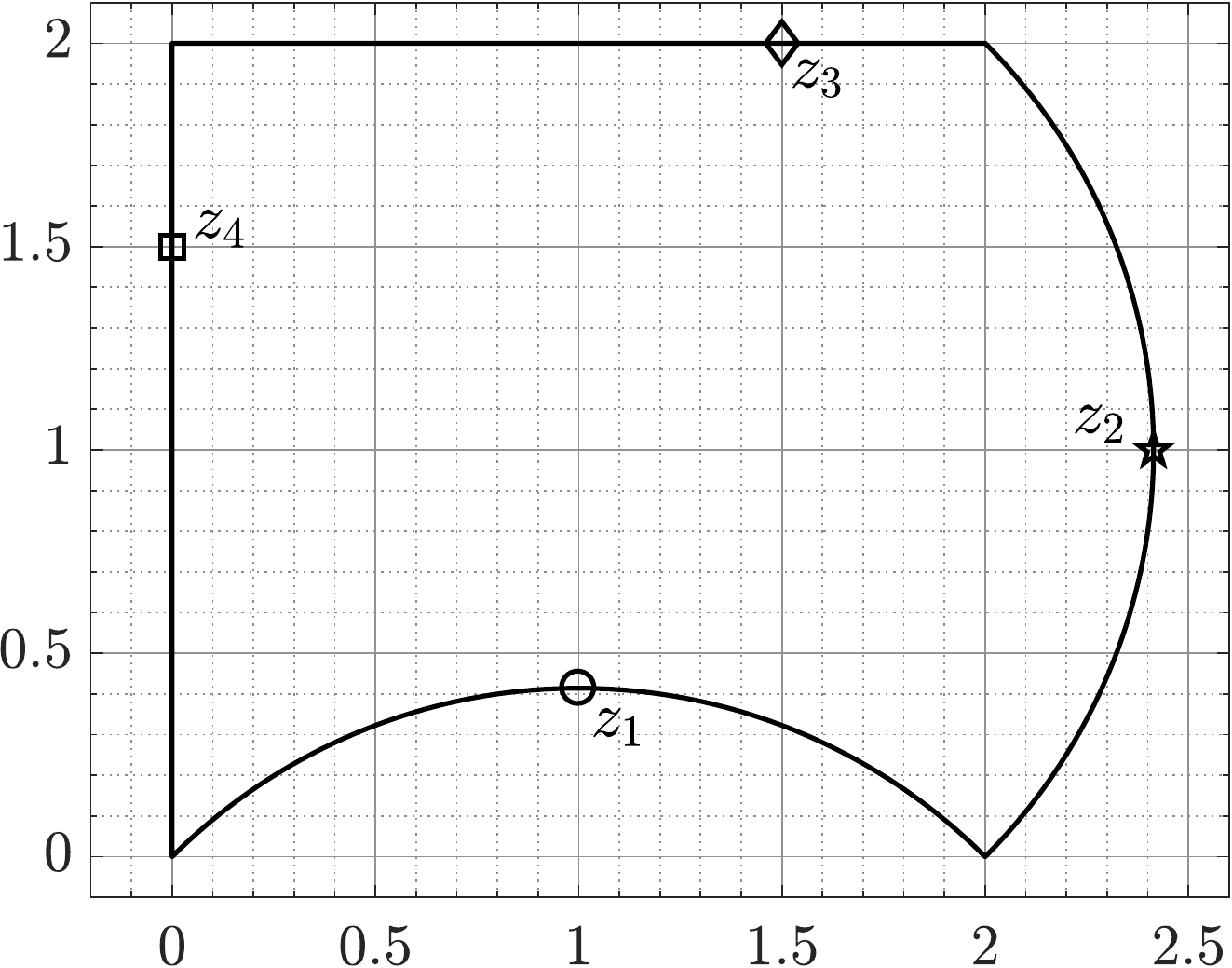}}
\hfill
\scalebox{0.55}{\includegraphics[trim=0cm 0cm 0cm 0cm,clip]{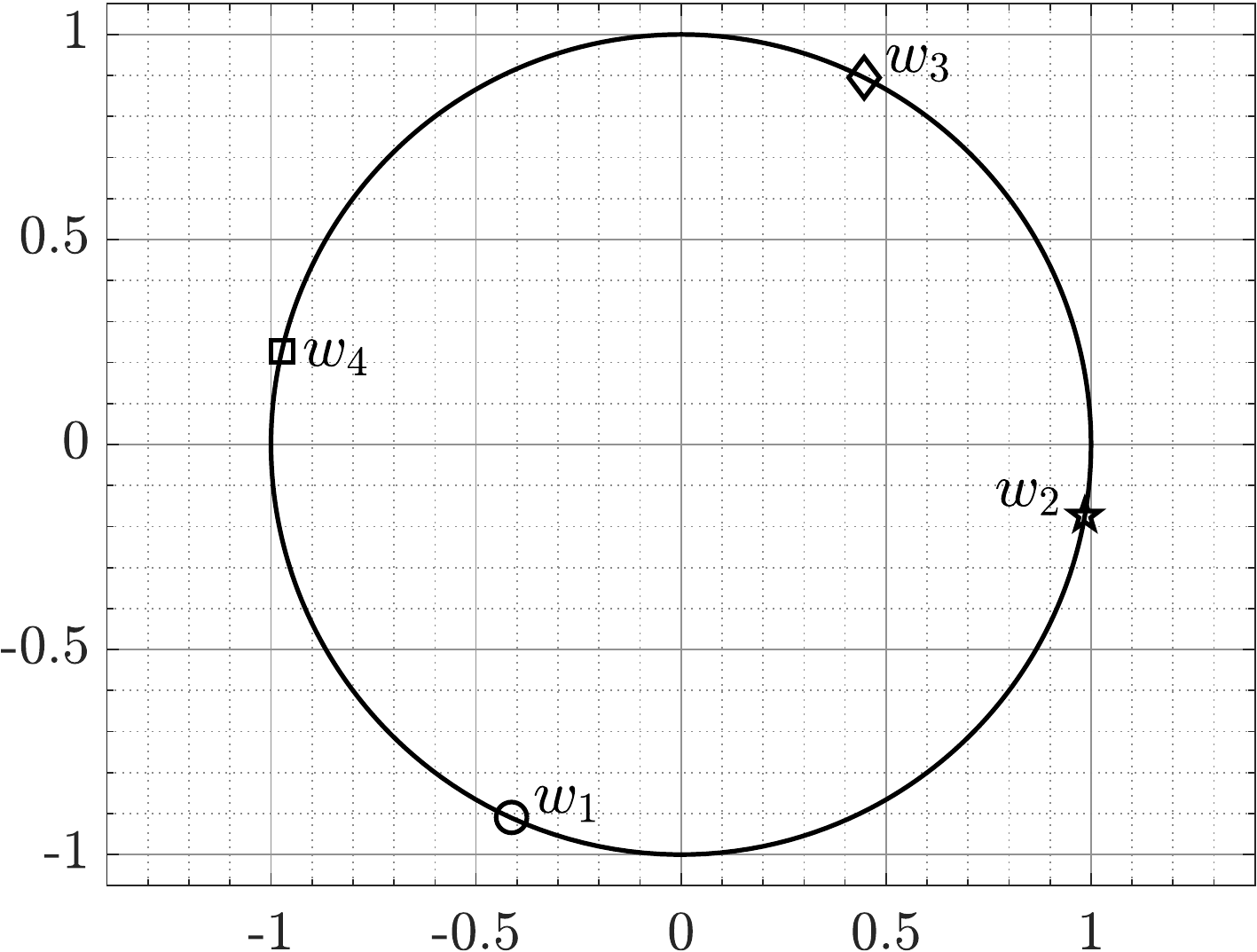}}
}
\caption{The quadrilateral $(D;z_1,z_2,z_3,z_4)$ (left) and the quadrilateral $(\DD;w_1,w_2,w_3,w_4)$ (right).}
  \label{fig:1}
\end{figure}

\end{mysubsection}

\begin{mysubsection}{\bf Exterior modulus of quadrilaterals.}
Let $D$ be a bounded simply connected polygonal domain and let $D^-=\overline{\CC}\backslash D$ be its complement with respect to the extended complex plane $\overline{\CC}=\CC\cup\{\infty\}$. 
We consider the quadrilateral $(D;z_1,z_2,z_3,z_4)$ where $z_1,z_2,z_3,z_4$ are four distinguished points on $\Gamma=\partial D$ with clockwise orientation. 
The exterior modulus of a quadrilateral $(D;z_1,z_2,z_3,z_4)$ equals the modulus of the family of all curves in $D^-$ joining the opposite boundary arcs $(z_2,z_3)$ and $(z_4,z_1)$.

As in the bounded case, the unbounded domain $D^-$ can be mapped using the conformal mapping $w=\Phi(z)$ described in \S~\ref{sec:map} onto the unit disk $\DD$ so that the four points $z_1,z_2,z_3,z_4$ on $\partial D$ (in clockwise orientation) are mapped onto four points $w_1,w_2,w_3,w_4$ on the unit circle $\partial\DD$ (in counterclockwise orientation). Assume that $z_k=\eta(\hat t_k)$ where $\hat t_k\in[0,2\pi]$, then~\eqref{eq:S} implies that $w_k=e^{\i S(\hat t_k)}$, $k=1,2,3,4$. By the conformal invariance of the modulus, the exterior modulus of the quadrilateral $(D;z_1,z_2,z_3,z_4)$ is equal to ${\rm mod}(\DD;w_1,w_2,w_3,w_4)$ which can be computed using the exact formula~\eqref{eq:mod-disk}.

\end{mysubsection}

\section{Examples: modulus of quadrilaterals}

\begin{nonsec}{\bf Half-disk.} \label{half} Let 
$$A=\{z \in {\mathbb C}: |z|<1, {\rm Im}z >0\},\quad -1<r<s<1, 0<\sigma<
\beta<\pi.$$ Then $(A; z_1,z_2,z_3,z_4)$ where $z_1=r$, $z_2=s$, $z_3=e^{\i\sigma}$, $z_4=e^{\i\beta}$, is a quadrilateral and its modulus can be computed by
means of elementary conformal mappings and it is
\begin{equation}\label{halfmod}
{\rm mod}(A; z_1,z_2,z_3,z_4)=\frac{\pi}{2} \frac{1}{\mu(1/\sqrt{u})}\,,
\end{equation}
where
\begin{equation}\label{halfmod2}
 u=|f_1(\exp(\i\beta))^2, f_1(r)^2,f_1(s)^2, f_1(\exp(\i\sigma))^2 |\,,\quad f_1(z) = \frac{1+z}{1-z}\,.
\end{equation}

\begin{figure}[H] %
\centerline{
\scalebox{0.7}{\includegraphics[trim=0cm 0cm 0cm 0cm,clip]{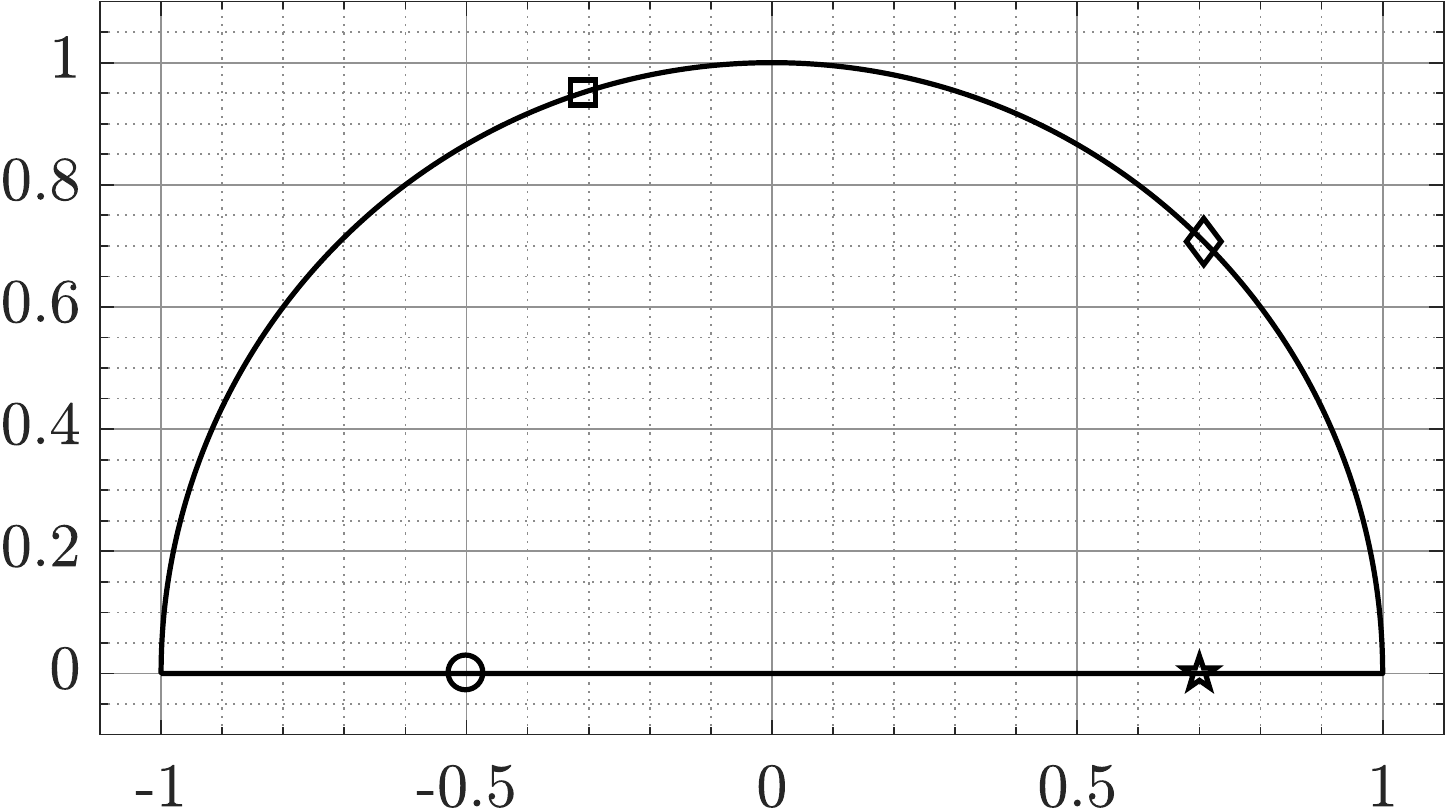}}
}
\caption{Half disk quadrilateral with $z_1=-0.5$, $z_2=0.7$, $z_3=e^{\pi\i/4}$, and $z_4=e^{3\pi\i/5}$.}
\label{fig:halfdisk}
\end{figure}

The approximate numerical values of the modulus ${\rm mod}(A; z_1,z_2,z_3,z_4)$ are computed by the above described method with $n=2^{13}$ and the exact values are presented in Table~\ref{tab:half}. The relative error in the approximated values are also presented in Table~\ref{tab:half} which is of order $10^{-14}$.

\begin{table}[H]
 \begin{tabular}{|l|l|l|l|l|l|c|} 
 \hline
   $r$ & $s$  & $m$  & $n$  & Exact Value   & Approx. Value  & Relative Error \\ \hline
 -0.8  & 0.2  & 1/8  & 1/3  & 1.12275580085474 	& 1.12275580085477 	& $2.71\times10^{-14}$ \\ \hline
 -0.5  & 0.3  & 1/6  & 2/5  & 0.96809243696619 	& 0.96809243696621 	& $2.10\times10^{-14}$\\ \hline
 -0.2  & 0.5  & 1/5  & 1/2  & 0.79872083257913 	& 0.79872083257913 	& $1.14\times10^{-14}$ \\ \hline
  0.2  & 0.6  & 1/4  & 3/5  & 0.95886428362598 	& 0.95886428362597 	& $1.83\times10^{-14}$\\ \hline
  0.2  & 0.8  & 1/4  & 4/5  & 0.83635871682559 	& 0.83635871682557 	& $2.83\times10^{-14}$ \\ \hline
\end{tabular}
 \caption{The values of the modulus of the quadrilateral for several values of $r,s,\sigma=m\pi, \beta=n\pi$.}\label{tab:half} 
\end{table}

\end{nonsec}

\begin{nonsec}{\bf Trapezoid.}
In this example, we consider the trapezoid $T$ with the vertices $z_1=0$, $z_2=1$, $z_3=1+\i L$, $z_4=\i(L-1)$ (see Figure\ref{fig:Trapezoid} (left)). The exact value of the modulus ${\rm mod}(T; z_1,z_2,z_3,z_4)$ is given for $L>1$ by~\cite[p.~82]{ps10}
\begin{equation}\label{eqn:mod-trap}
{\rm mod}(T; z_1,z_2,z_3,z_4)=\frac{\pi}{2\mu(\kappa)},  
\end{equation}
where
\[
\kappa=\frac{1-2\lambda\lambda'}{1+2\lambda\lambda'}, \quad
\lambda=\mu^{-1}\left(\frac{\pi}{2(2L-1)}\right),\quad \lambda'=\sqrt{1-\lambda^2}.
\]

The above method  with $n=2^{13}$ is used to compute approximate values of the modulus ${\rm mod}(T; z_1,z_2,z_3,z_4)$ for several values of $L\in (1,5]$. The relative error in the computed values is presented in Figure~\ref{fig:Trapezoid}. 
The values of the function $\mu$ and its inverse are computed as described in~\cite{nv}.
Figure~\ref{fig:Trapezoid} presents also the relative error in the approximate values obtained using the SC Toolbox~\cite{Dri}. As we see from the figures, for $L>2.5$, the relative errors of the two methods are almost identical. For small $L$, the results obtained with SC Toolbox are better than the results obtained by the proposed method.

\begin{figure}[H] %
\centerline{
\scalebox{0.7}{\includegraphics[trim=0cm 0cm 0cm 0cm,clip]{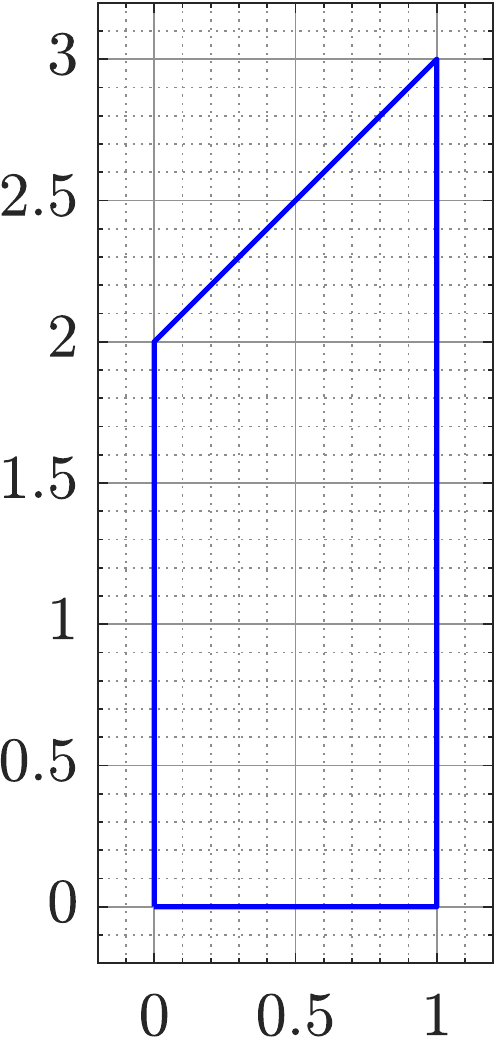}}
\hfill
\scalebox{0.7}{\includegraphics[trim=0cm 0cm 0cm 0cm,clip]{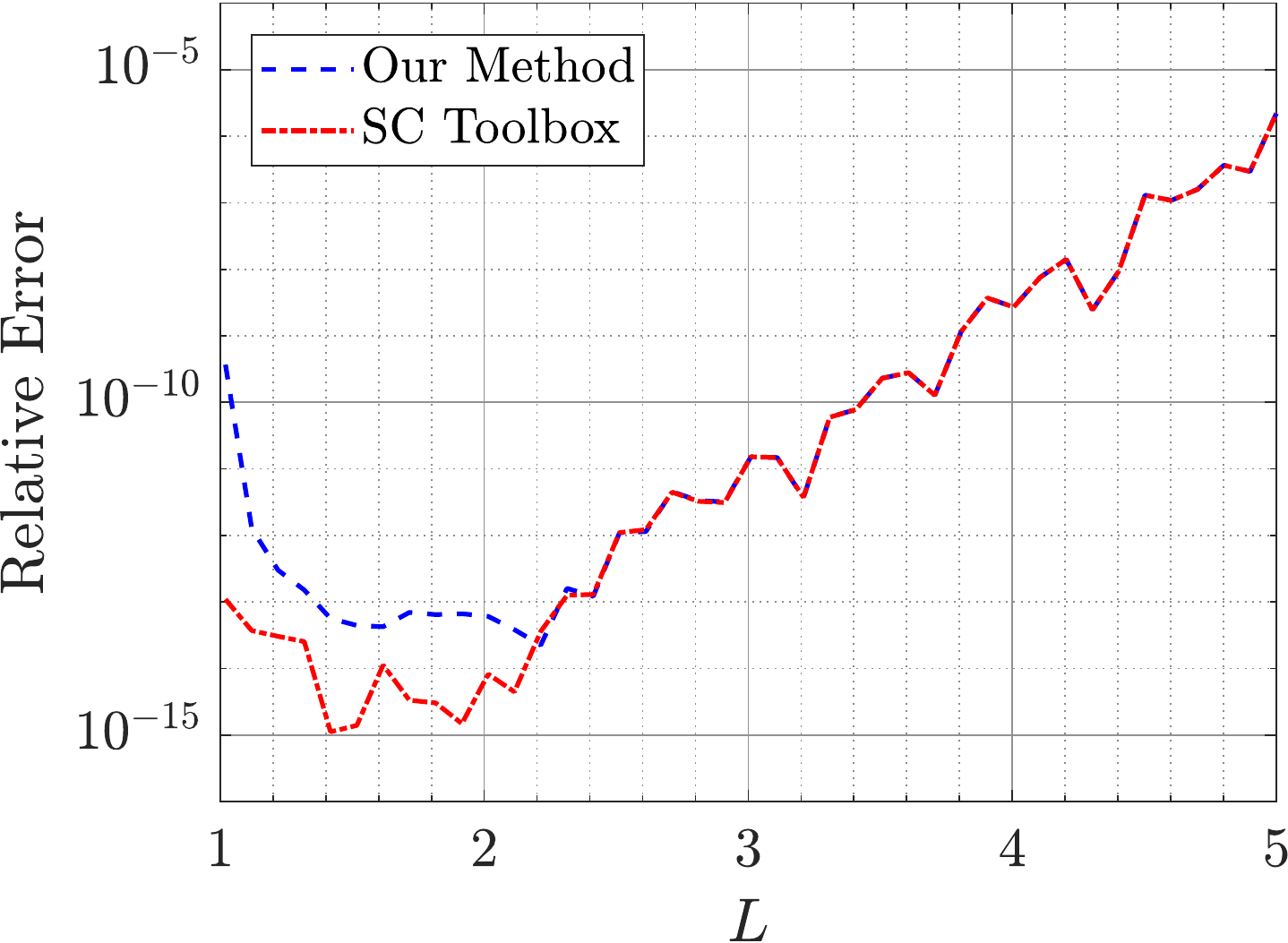}}
}
\caption{The Trapezoid $T$ for $L=3$ (left) and the relative error in the computed approximate values of the modulus vs $L$ for the trapezoid (right).}
\label{fig:Trapezoid}
\end{figure}

\end{nonsec}

\begin{nonsec}{\bf Symmetric trapezoid.}
Let $a,b>0$ and consider the polygon $P$ with the vertices $z_1=-a+\i$, $z_2=-b$, $z_3=b$, and $z_4=a+\i$. Then, by symmetry, 
\[
{\rm mod}(P; z_1,z_2,z_3,z_4)=2{\rm mod}(P_+; \i,0,b,a+\i)
\]
where $P_+$ is the polygon with vertices $\i$, $0$, $b$, and $a+\i$ (see Figure~\ref{fig:Trap-s1} (left)). If $b=a+1$, then by~\eqref{eqn:mod-trap}, the exact value of the modulus is given by 
\[
{\rm mod}(R; z_1,z_2,z_3,z_4)=\frac{\pi}{\mu(\kappa)}, \quad
\kappa=\frac{1-2\lambda\lambda'}{1+2\lambda\lambda'}, \quad
\lambda=\mu^{-1}\left(\frac{\pi}{2(2a+1)}\right),\quad \lambda'=\sqrt{1-\lambda^2}.
\]

For $b=a+1$, the proposed method is used with $n=2^{13}$ to compute approximate values of ${\rm mod}(P_+; \i,0,z_3,z_4)$ for several values of $a$. The relative error in the computed values is presented in Figure~\ref{fig:Trap-s1} (right). 
For other values of $b$, let the real function $u(a,b)$ be defined for $(a,b)\in[1,3]\times[1,3]$ by
\begin{equation}\label{eq:trap-u}
u(a,b) = 2{\rm mod}(P_+; \i,0,b,a+\i).
\end{equation}
The values of the function $u(a,b)$ are computed for several values of $a,b$ using the above proposed method and using the SC Toolbox. The contour lines of the values of the function $u$ for the proposed method are presented in Figure~\ref{fig:Trap-s2} (left). In Figure~\ref{fig:Trap-s2} (right), we present the absolute value of the difference between the values of the function $u$ obtained by the above method and by the SC Toolbox.

\begin{figure}[H] %
\centerline{
\scalebox{0.6}{\includegraphics[trim=0cm 0cm 0cm 0cm,clip]{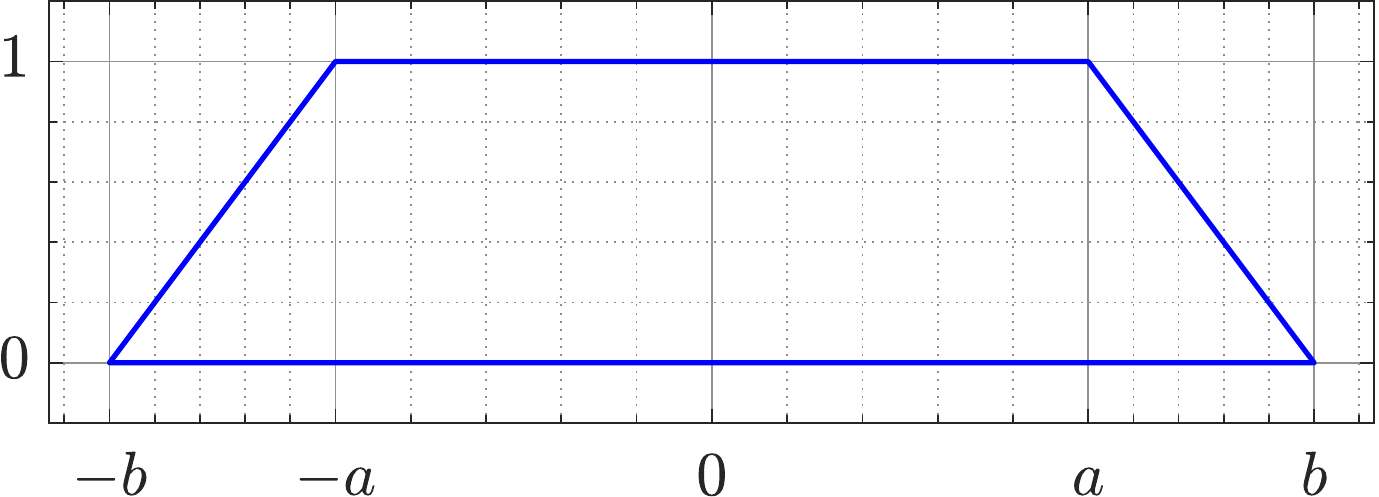}}
\hfill
\scalebox{0.6}{\includegraphics[trim=0cm 0cm 0cm 0cm,clip]{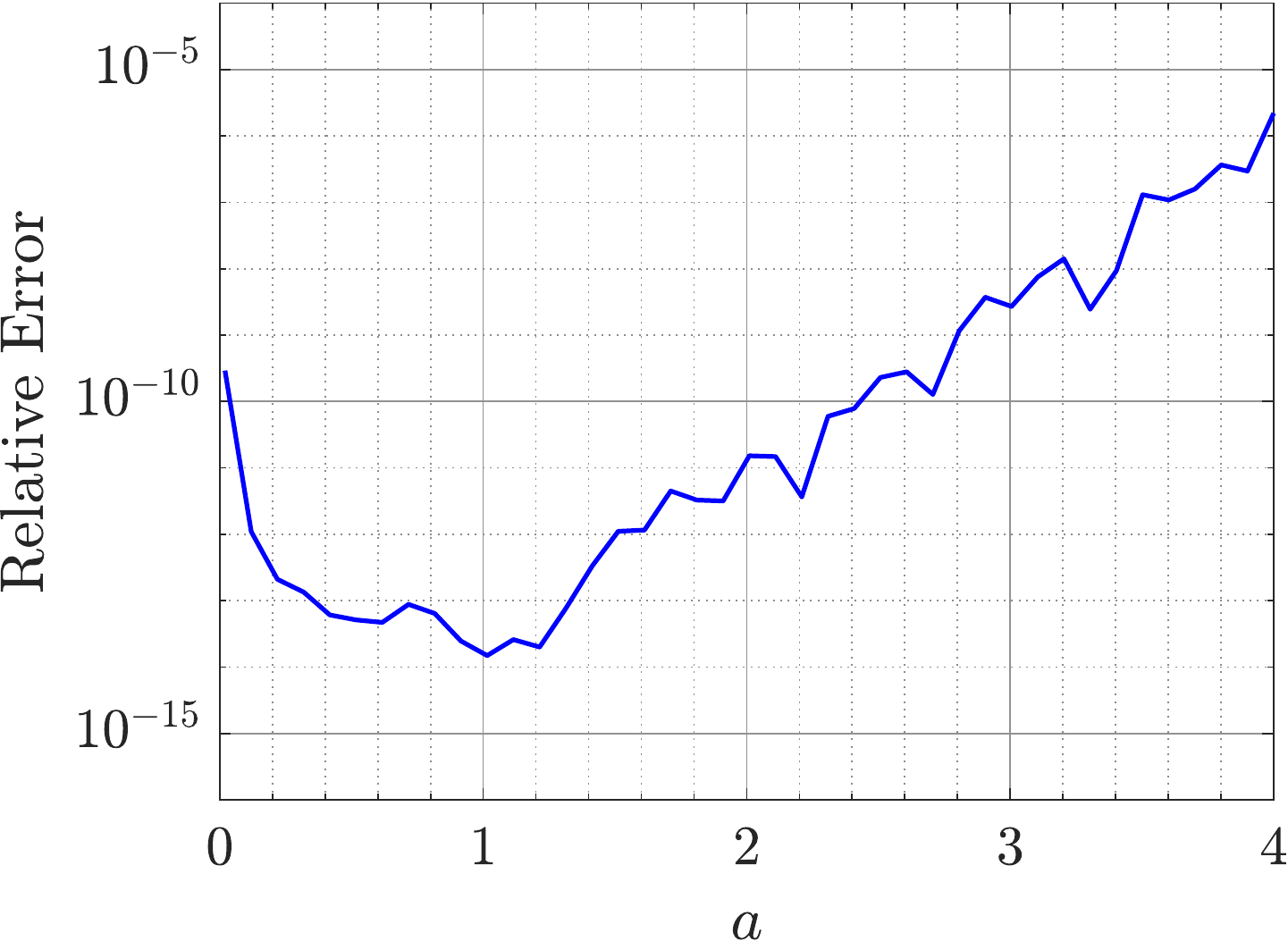}}
}
\caption{The symmetric trapezoid $P$ (left) and the relative error in the computed approximate values of the modulus vs $a$ (right).}
\label{fig:Trap-s1}
\end{figure}

\begin{figure}[H] %
\centerline{
\scalebox{0.6}{\includegraphics[trim=0cm 0cm 0cm 0cm,clip]{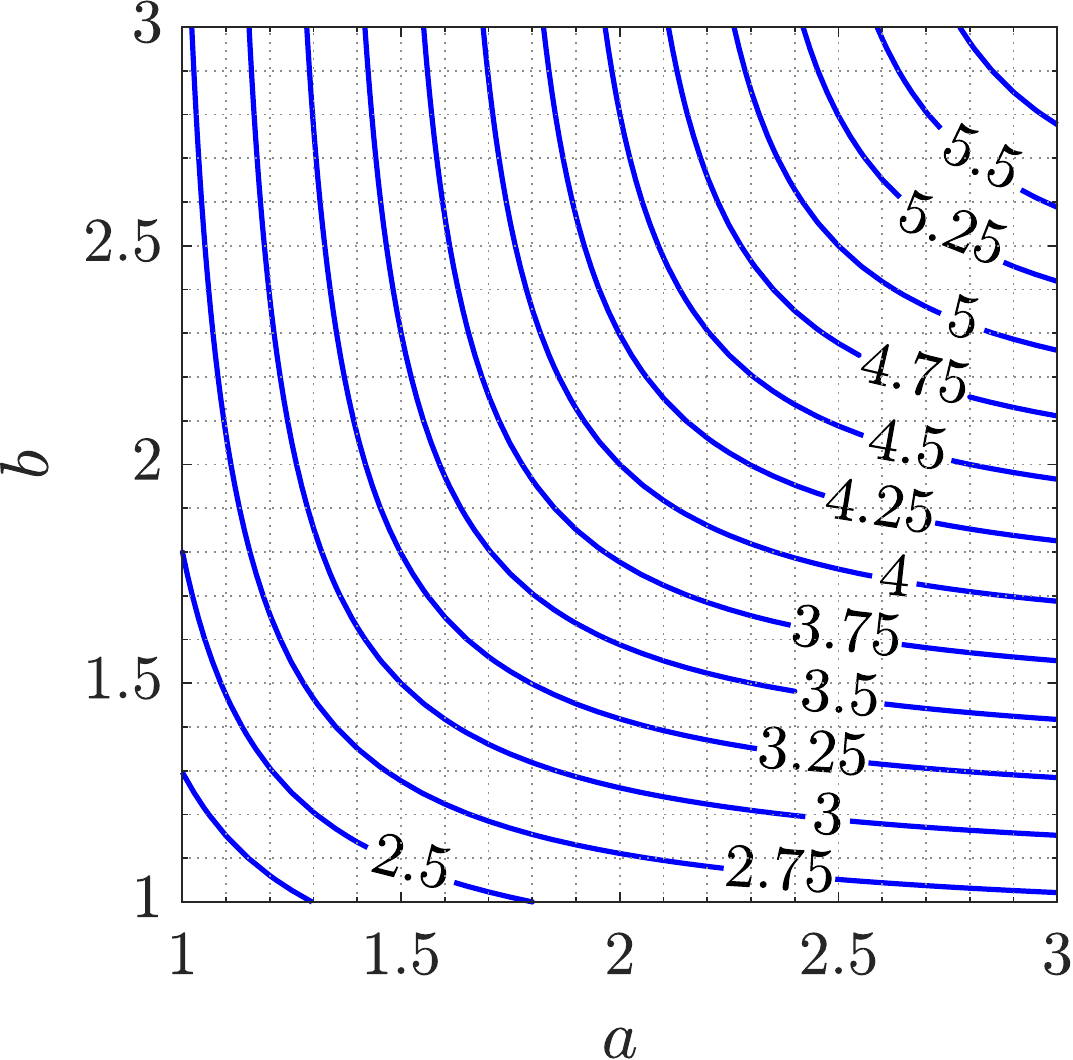}}
\hfill
\scalebox{0.6}{\includegraphics[trim=0cm 0cm 0cm 0cm,clip]{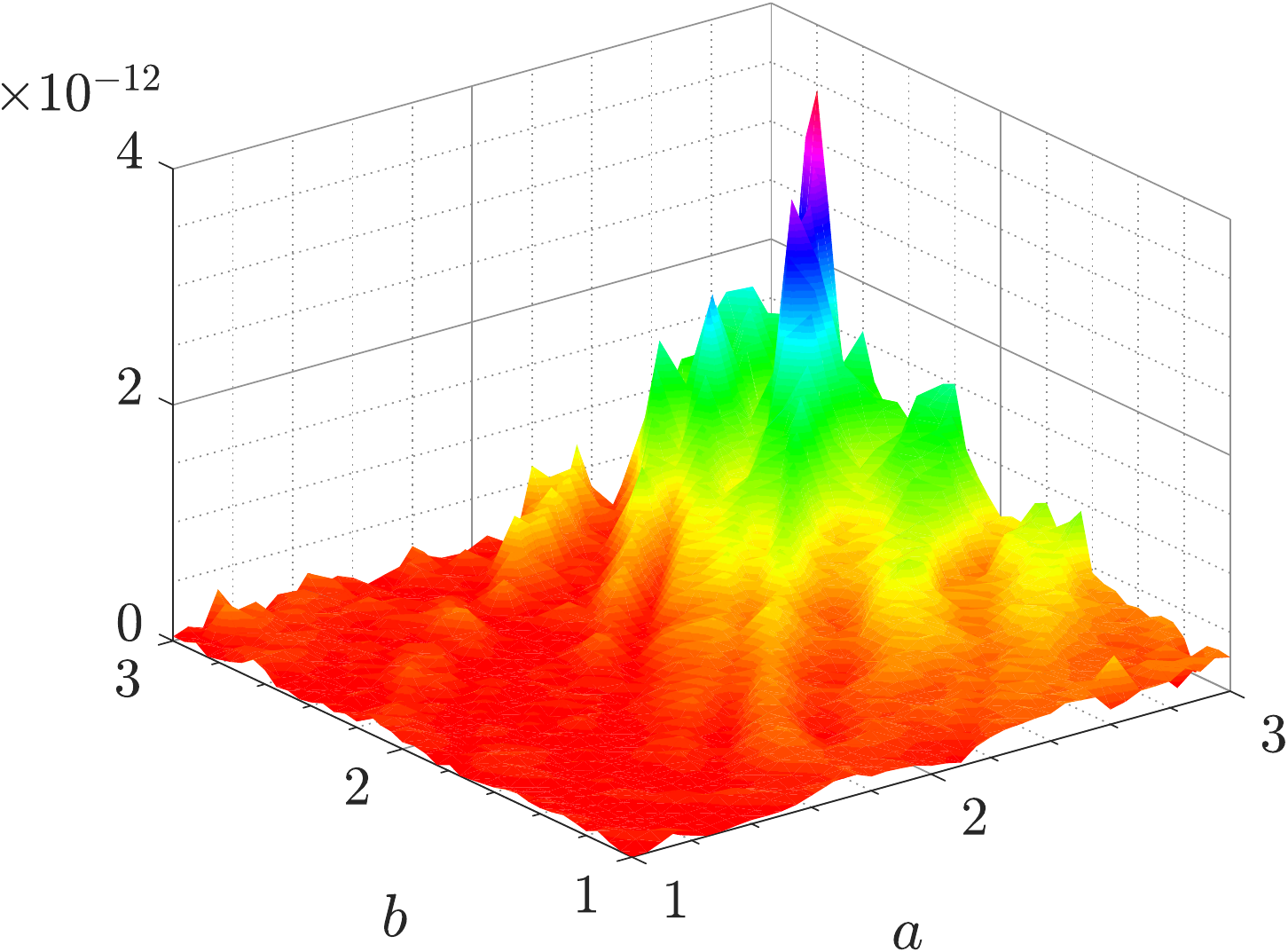}}
}
\caption{The contour lines of the function $u(a,b)$ (left) the absolute value of the difference between the values of the function $u$ obtained by the presented method and by the SC Toolbox.}
\label{fig:Trap-s2}
\end{figure}

\end{nonsec}

\begin{nonsec}{\bf Trapezoid with a small curvature at the left side.}
In this example, we consider the trapezoid $T_\epsilon$ with the vertices $z_1=0$, $z_2=1+\i\tan\pi\sigma$, $z_3=1+\i(L+\tan\pi\beta)$, $z_4=\i L$, and with a small curvature at the left side parametrized the angle $\pi\epsilon$ (see Figure~\ref{fig:Trap1} (left)). This trapezoid has been considered in~\cite[p.~14]{acnc} but no numerical results were presented in~\cite{acnc}. 

In this paper, we use the above method with $n=2^{12}$ to compute approximate values of the modulus ${\rm mod}(T_\epsilon; z_1,z_2,z_3,z_4)$ for several values of $0.0001\le\epsilon\le0.25$. 
The values of the computed modulus is presented in Figure~\ref{fig:Trap1} (right) for $L=2$, $\sigma=\pi/8$ and $\beta=\pi/4$. 

\begin{figure}[H] %
\centerline{
\scalebox{0.7}{\includegraphics[trim=0cm 0cm 0cm 0cm,clip]{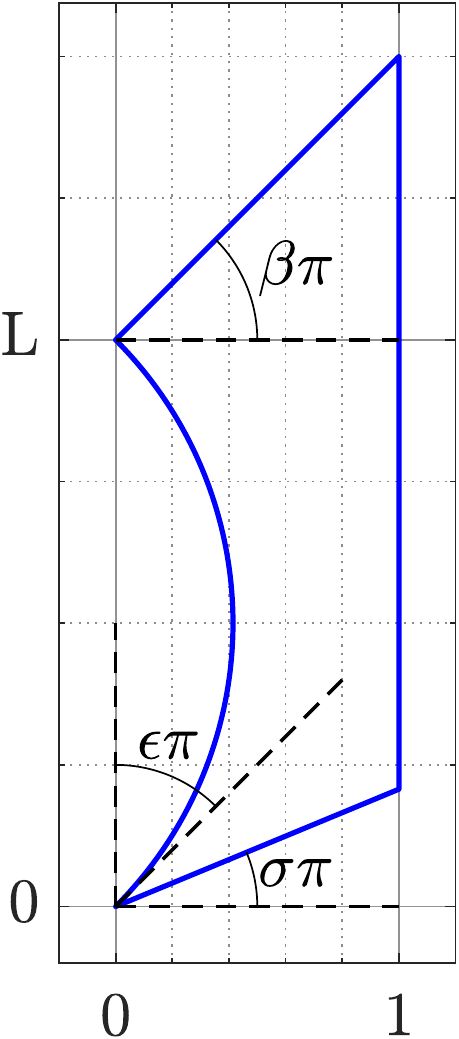}}
\hfill
\scalebox{0.7}{\includegraphics[trim=0cm 0cm 0cm 0cm,clip]{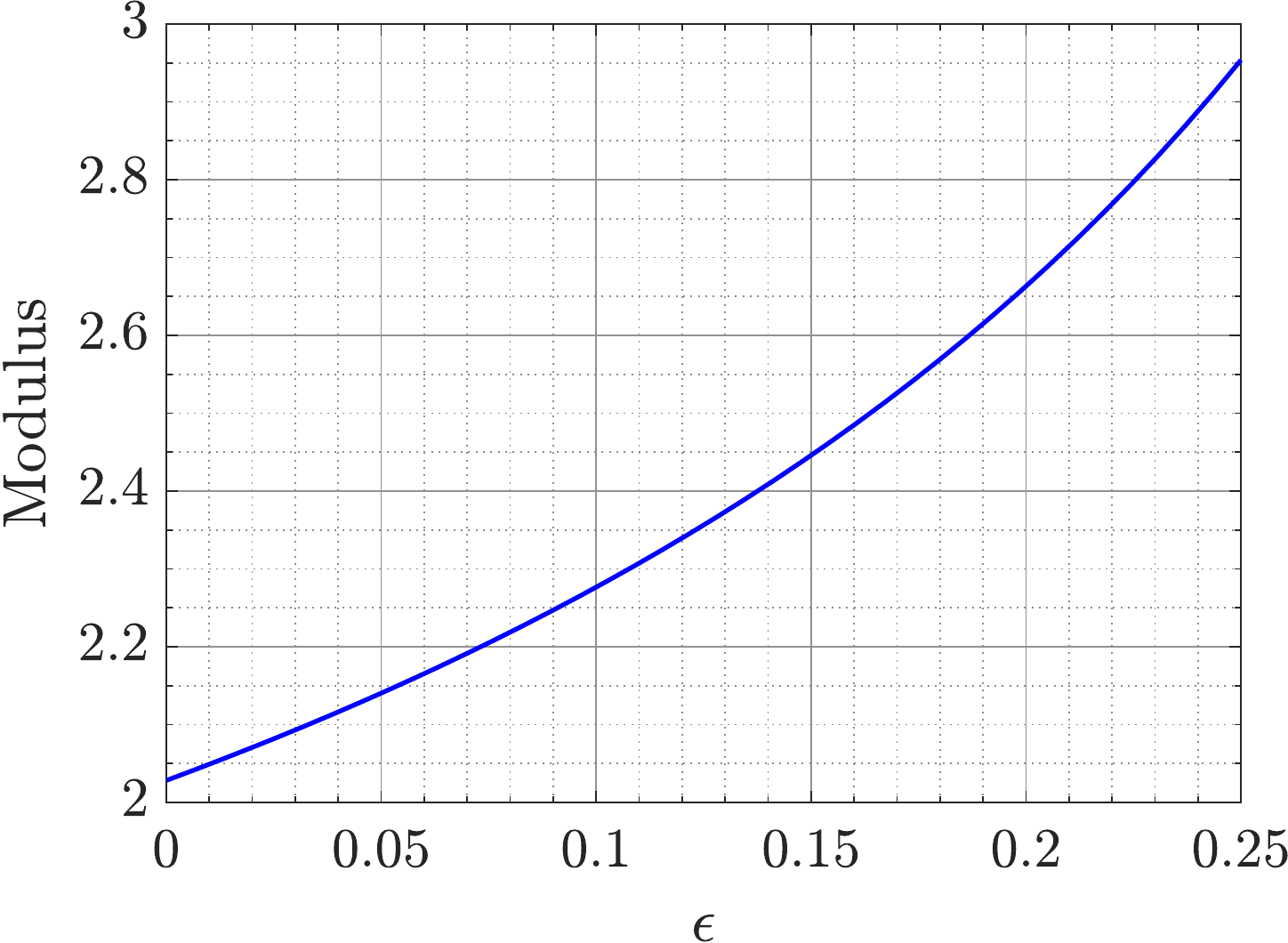}}
}
\caption{The Trapezoid $T_\epsilon$ (left) and the values of its modulus (right).}
\label{fig:Trap1}
\end{figure}

For $\epsilon=0$, the exact value of the capacity is given by~\cite[Eq.~(74)]{acnc}
\begin{equation}\label{eqn:crowdy}
{\rm mod}(T_\epsilon; z_1,z_2,z_3,z_4)=\frac{2\mu(\kappa)}{\pi}, \quad 
\kappa=\sqrt{t_0}
\end{equation}
where
\[
t_0=16e^{-\pi\hat{L}}\left(1-8(1+4\sigma\beta)e^{-\pi\hat{L}}+4(11+4\sigma^2+4\beta^2+128\sigma\beta+304\sigma^2\beta^2)e^{-2\pi\hat{L}}+O(e^{-3\pi\hat{L}})\right),
\]
\[
\hat{L}=L+\frac{1}{\pi}\left(\Psi(0.5+\beta)+\Psi(0.5-\sigma)-2\Psi(0.5)\right),
\]
and $\Psi(z)=\Gamma'(z)/\Gamma(z)$ is the digamma function. Here, we approximate $t_0$ with
\[
t_0\approx 16e^{-\pi\hat{L}}\left(1-8(1+4\sigma\beta)e^{-\pi\hat{L}}+4(11+4\sigma^2+4\beta^2+128\sigma\beta+304\sigma^2\beta^2)e^{-2\pi\hat{L}}\right).
\]
Figure~\ref{fig:Trap2} presents the absolute values of the difference between the approximate values of the modulus obtained with our method and with the Formula~\eqref{eqn:crowdy} for several values of $0\le\sigma\le0.25$ and $\epsilon=0$, $L=2$, $\beta=0.25$. It follows from the definition of $\Hat L$ that the values of $\hat L$ decrease as $\sigma$ increases. Thus, it is expected that the difference between the approximate values of the modulus increase as $\sigma$ increases since $t_0$ is approximated to within $O(e^{-3\pi\hat{L}})$.

\begin{figure}[H] %
\centerline{
\scalebox{0.7}{\includegraphics[trim=0cm 0cm 0cm 0cm,clip]{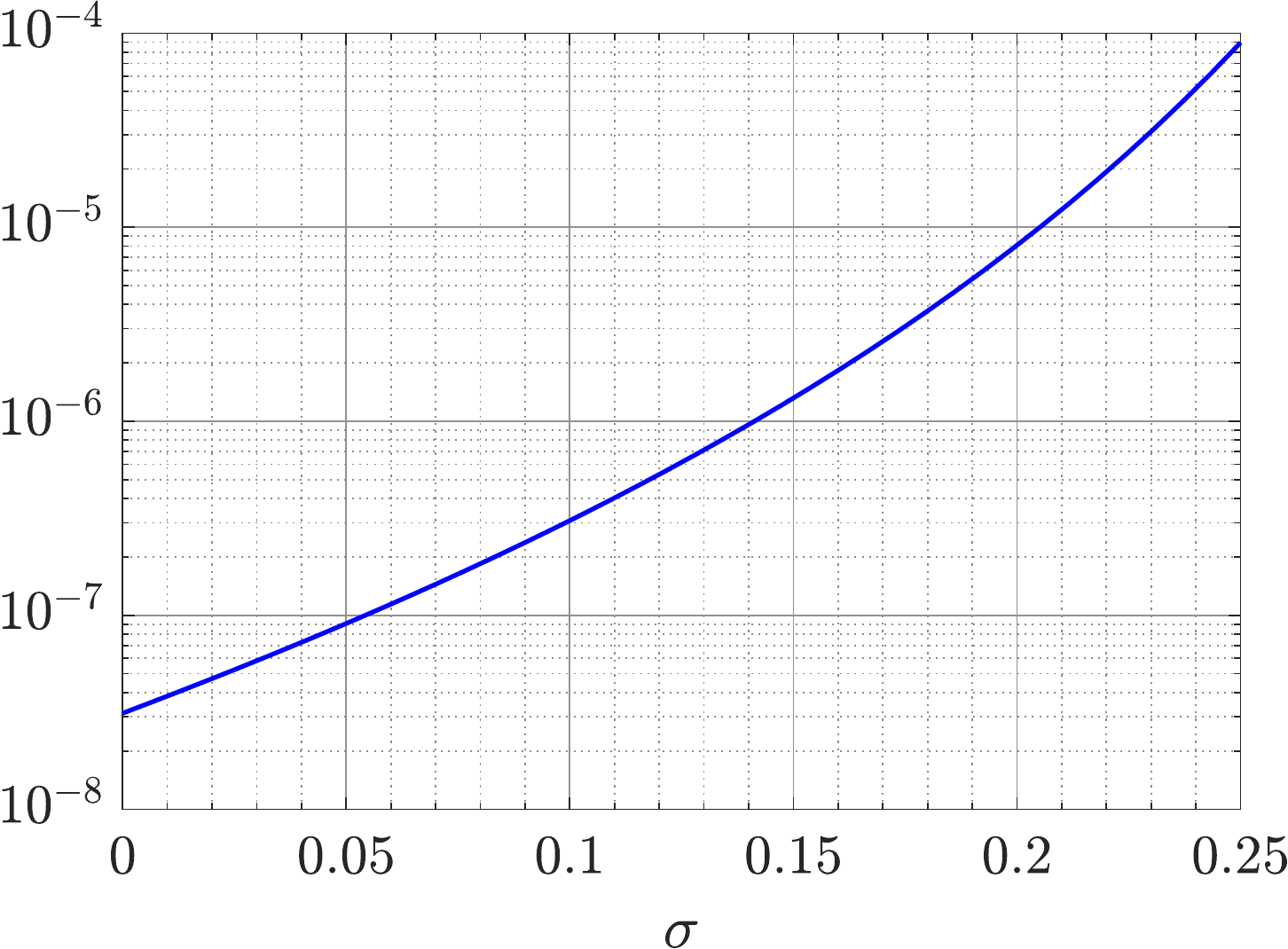}}
}
\caption{The absolute values of the difference between the approximate values of the modulus obtained with our method and with the Formula~\eqref{eqn:crowdy}.}
\label{fig:Trap2}
\end{figure}

\end{nonsec}

\begin{nonsec}{\bf $L$-shaped domains: polygonal boundary.}\label{exm:Lp}
Consider the simply connected domain interior to the polygon with the vertices $v_1=-1+3\i$, $v_2=-1+\i$, $v_3=-1-\i$, $v_4=1-\i$, $v_5=3-\i$, $v_6=3+\i$, $v_7=1+\i$, and $v_8=1+3\i$ (see~\cite[p.~44]{ps10}). We consider here $v_2$ and $v_4$ as vertices. 

\begin{figure}[H] %
\centerline{
\scalebox{0.6}{\includegraphics[trim=0cm 0cm 0cm 0cm,clip]{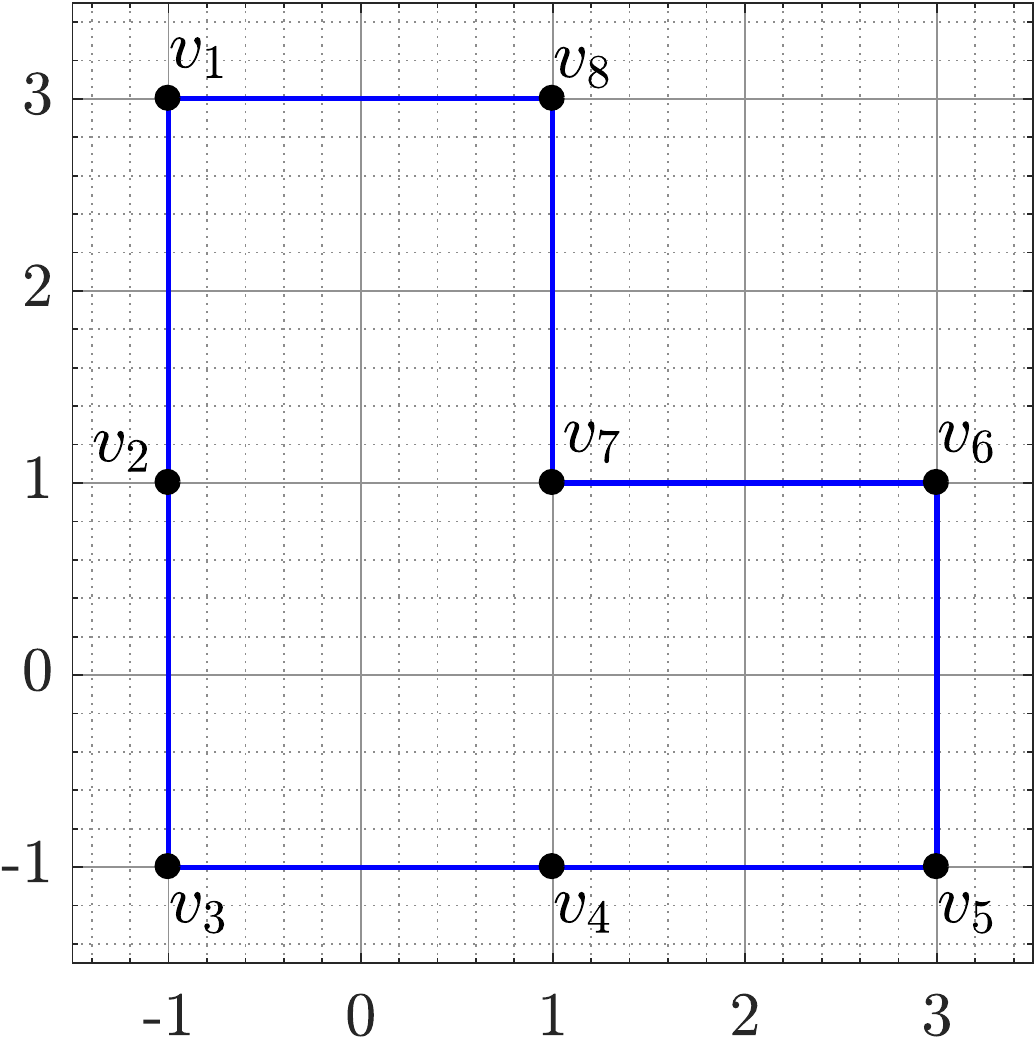}}
\hfill
\scalebox{0.6}{\includegraphics[trim=0cm 0cm 0cm 0cm,clip]{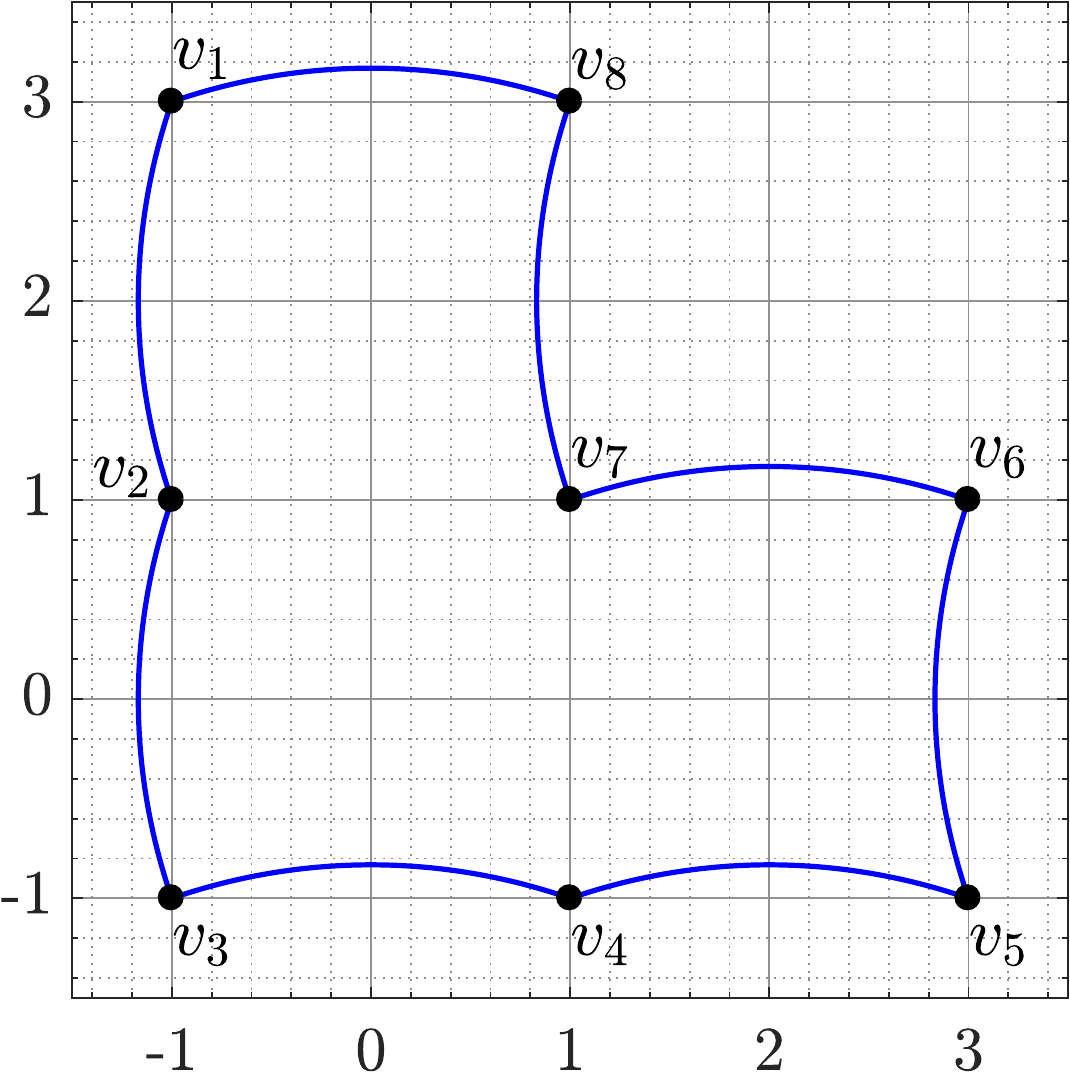}}
}
\caption{Left: $L$-shaped quadrilateral. Right: The same with
circular arc boundary curves.}
\label{fig:letterL}
\end{figure}

The $L$-shaped domain with four vertices (in counterclockwise orientation) is a quadrilateral. There are $280$ possible choices of such quadrilaterals. The proposed method is used with $n=2^{13}$ to compute the modulus for these $280$ quadrilaterals and the reciprocal error based on the identity~\eqref{reciprel}, i.e., the \emph{reciprocal error} in Tables~\ref{tab:L} and~\ref{tab:Larc} is defined by
\begin{equation}\label{eq:rec-err}
\left|1-{\rm mod}(L;z_1, z_2, z_3, z_4)\,{\rm mod}(L;z_2, z_3, z_4,z_1)\right|.
\end{equation}
The results are presented in Figure~\ref{fig:L-reci} (left). 

Extensive numerical  tests \cite{HRV} related to capacity computation show that the reciprocal error~\eqref{eq:rec-err} agrees with several other error estimates. However, for the current proposed method, the reciprocal error is not significant because the method is based on mapping the domain $D$ and the four points $z_1,z_2,z_3,z_4$ on its boundary to the unit disk $\DD$ with four points $w_1,w_2,w_3,w_4$ on the unit circle. 
Then ${\rm mod}(\DD,w_1,w_2,w_3,w_4)$ is computed using the exact formula. 
Thus the reciprocal error in ${\rm mod}(D,z_1,z_2,z_3,z_4)$ is the same as the reciprocal error in  ${\rm mod}(\DD,w_1,w_2,w_3,w_4)$. Thus, the reciprocal error in our method measures only the error in the numerical computation of the special function ``$\mu$'' in the exact formula.

For this example, the exact values of the modulus of the $L$-shaped quadrilateral for several choices of vertices are given in~\cite{Ga}. Table~\ref{tab:L} presents these exact values as well as the approximate values obtained using the proposed method and the relative error in the approximate values. 
In this example, as well as in the next example, the auxiliary point $\alpha$ in~\eqref{eq:A} need to be chosen carefully to ensure the convergence of the method. In our numerical computation we choose $\alpha$ inside the domain $L$, sufficiently far from the boundary, and close to the arithmetic mean of the vertices $z_1$, $z_2$, $z_3$, $z_4$ of the quadrilateral $(L;z_1, z_2, z_3, z_4)$.

\begin{figure}[H] %
\centerline{
\scalebox{0.55}{\includegraphics[trim=0cm 0cm 0cm 0cm,clip]{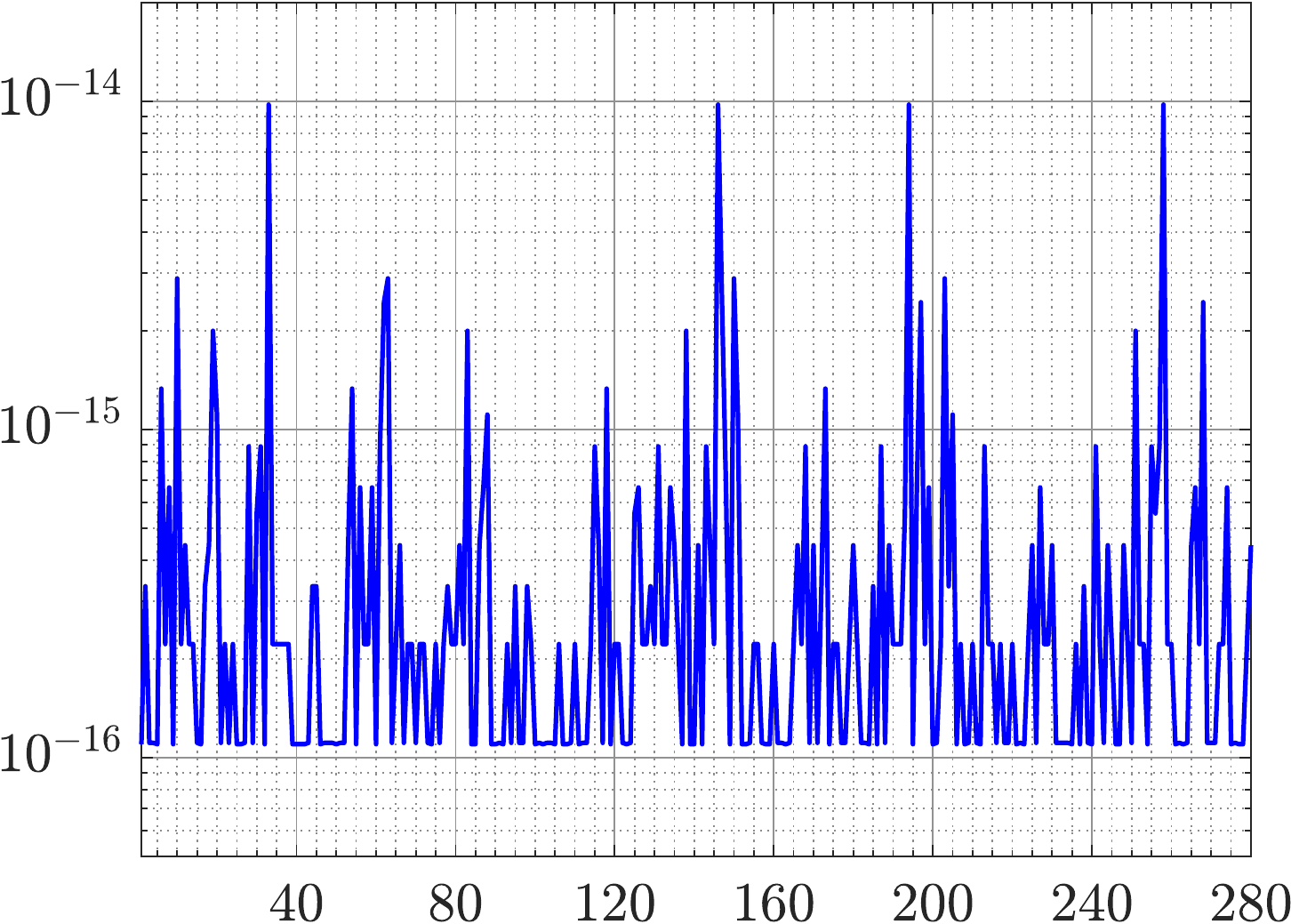}}
\hfill
\scalebox{0.55}{\includegraphics[trim=0cm 0cm 0cm 0cm,clip]{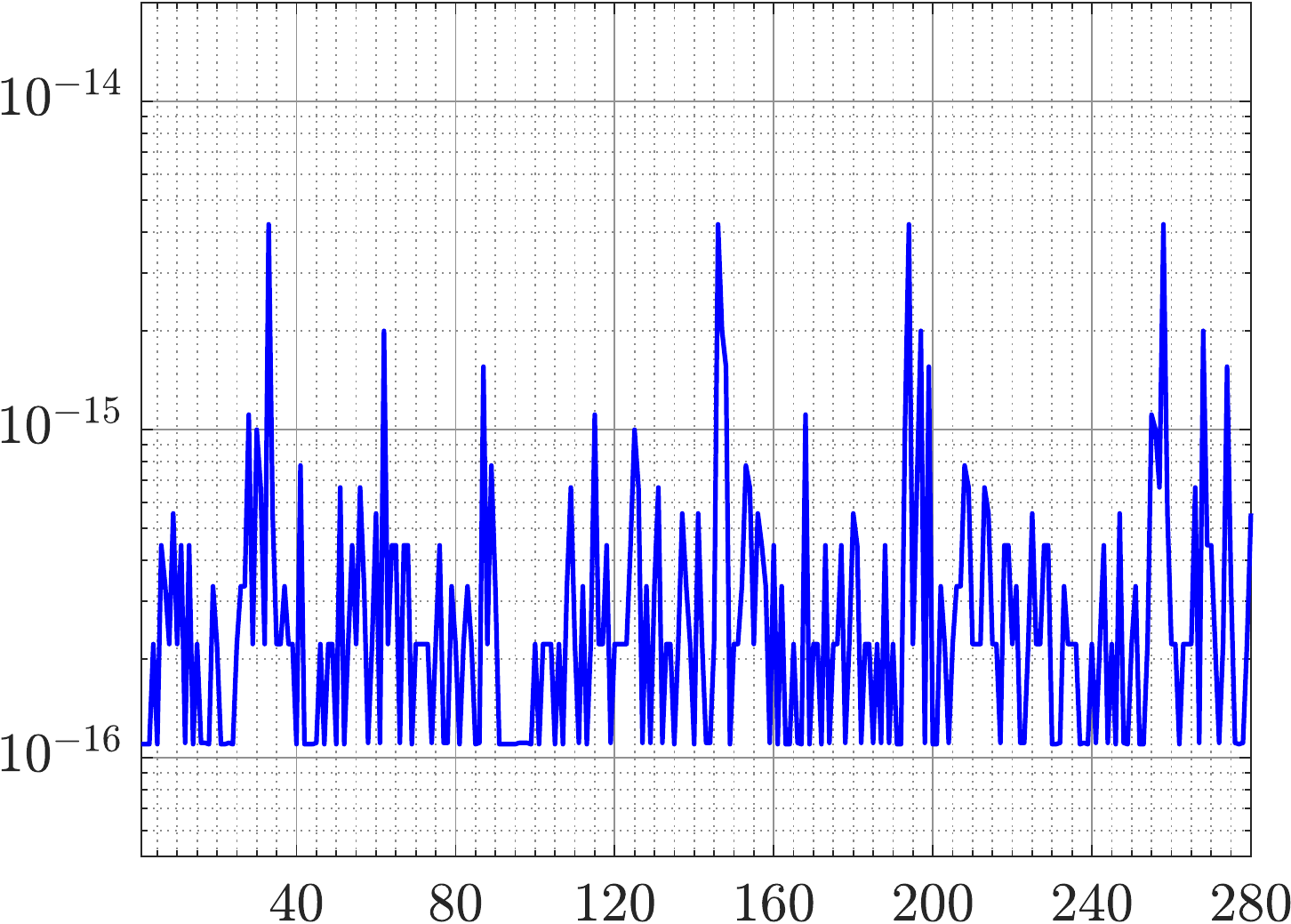}}
}
\caption{Left: The reciprocal error for the $L$-shaped quadrilateral. Right: The same with circular arc boundary curves.}
\label{fig:L-reci}
\end{figure}

\begin{table}[H]
  \caption{The values of the modulus of the $L$-shaped quadrilateral for several choices of vertices.}\label{tab:L}  
	\vspace{0.2cm}
\begin{tabular}{|l|l|l|c|c|} 
 \hline
$z_1, z_2, z_3, z_4$  & ${\rm mod}(L;z_1, z_2, z_3, z_4)$  & Exact Modulus & Relative Error & Reciprocal Error\\ \hline
$v_1,v_3,v_5,v_6$& $1.73205080757599$ & $1.73205080756888$ & $4.11\times10^{-12}$ & $1.11\times10^{-15}$ \\ \hline
$v_5,v_6,v_7,v_8$& $1.73205080756327$ & $1.73205080756888$ & $3.24\times10^{-12}$ & $4.44\times10^{-16}$ \\ \hline
$v_1,v_3,v_5,v_7$& $0.99999999999876$ & $1.00000000000000$ & $1.24\times10^{-12}$ & $1.11\times10^{-16}$ \\ \hline
$v_8,v_4,v_6,v_7$& $0.78170096135130$ & $0.78170096134806$ & $4.15\times10^{-12}$ & $2.22\times10^{-16}$ \\ \hline
$v_8,v_1,v_3,v_6$& $1.70916888656242$ & $1.70916888655749$ & $2.88\times10^{-12}$ & $2.22\times10^{-16}$ \\ \hline
$v_8,v_1,v_5,v_6$& $2.55852314235398$ & $2.55852314234188$ & $4.73\times10^{-12}$ & $9.77\times10^{-16}$ \\ \hline
$v_8,v_2,v_4,v_6$& $1.56340192270204$ & $1.56340192269611$ & $3.79\times10^{-12}$ & $6.66\times10^{-16}$ \\ \hline
\end{tabular}
   \end{table}

\begin{table}[H]
  \caption{The values of the modulus of the $L$-shaped (circular arc) quadrilateral for several choices of vertices.}
\label{tab:Larc}  
 	\vspace{0.2cm}
  \begin{tabular}{|l|l|l|c|} 
 \hline
   $z_1, z_2, z_3, z_4$  & ${\rm mod}(L;z_1, z_2, z_3, z_4)$    \\ \hline
$v_1,v_3,v_5,v_6$& $1.74325313824307$    \\ \hline
$v_5,v_6,v_7,v_8$& $1.58841772274399$    \\ \hline
$v_1,v_3,v_5,v_7$& $1.10535075580239$    \\ \hline
$v_8,v_4,v_6,v_7$& $0.84849597438205$    \\ \hline
$v_8,v_1,v_3,v_6$& $1.72040886827649$    \\ \hline
$v_8,v_1,v_5,v_6$& $2.58094977005996$    \\ \hline
$v_8,v_2,v_4,v_6$& $1.57581541654770$    \\ \hline
\end{tabular}
\end{table}

\end{nonsec}

\begin{nonsec}{\bf $L$-shaped domains: circular arc polygonal boundary. }\label{exp:Lcr}
Now, we consider a circular arc polygon with the same vertices as in Example~\ref{exm:Lp} (see Figure~\ref{fig:letterL} (right)). This polygon is obtained by replacing each side-segments in the polygon in Figure~\ref{fig:letterL} (left) by a circular arc such that the angle between the segment and the tangent to the circular arc is $\epsilon$. We consider here $\epsilon=1/3$.
The proposed method is used with $n=2^{13}$ to compute the modulus for the $280$ possible choices of quadrilaterals and the reciprocal error based on the identity~\eqref{reciprel}. The obtained results are presented in Figure~\ref{fig:L-reci} (right). The values of the modulus of this $L$-shaped quadrilateral for several choices of vertices are given in Table~\ref{tab:Larc}.

\end{nonsec}

\begin{nonsec}{\bf Circular arc polygonal boundary. }
Consider the polygon with vertices $v_1=-2-2\i$, $v_2=0.4-2\i$, $v_3=1.4-2\i$, $v_4=2-2\i$, $v_5=2+0.8\i$, $v_6=-0.6+0.8\i$, and $v_7=-2-0.6\i$ (see~\cite[Fig. 6]{Tr}). This polygon consists of $5$ straight segments and $2$ circular arcs with centers $0.9-2\i$ and $-2+0.8\i$, respectively (see Figure~\ref{fig:Tre}).  

There are $140$ possible choices of four vertices $z_1, z_2, z_3, z_4$ to get a quadrilateral. The proposed method is used with $n=7\times 2^{10}$ to compute the modulus for these $140$ quadrilaterals and the approximate values of the modulus for some of these choices are presented in Table~\ref{tab:tre}.

\begin{figure}[H] %
\centerline{
\scalebox{0.6}{\includegraphics[trim=0cm 0cm 0cm 0cm,clip]{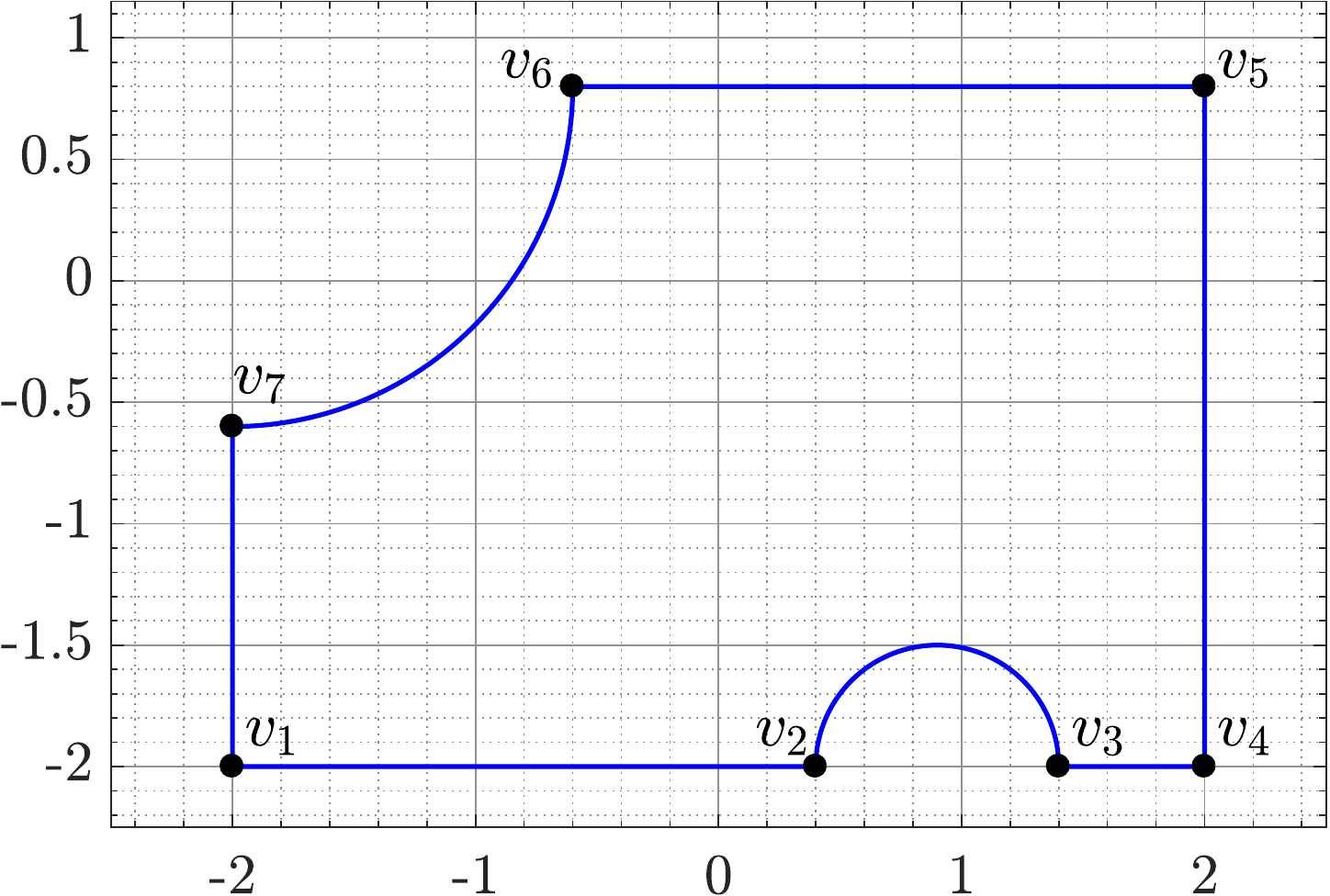}}
}
\caption{The polygon with $5$ straight segments and $2$ circular arcs.}
\label{fig:Tre}
\end{figure}

\begin{table}[H]
  \caption{The values of the modulus of the polygon with $5$ straight segments and $2$ circular arcs.}
\label{tab:tre}  
 	\vspace{0.2cm}
  \begin{tabular}{|l|l|c|} 
 \hline
$z_1,z_2,z_3,z_4$  & ${\rm mod}(L;z_1, z_2, z_3, z_4)$   \\ \hline
$v_1,v_2,v_3,v_4$& $2.45771442325834$   \\ \hline
$v_1,v_2,v_3,v_5$& $1.35593720891099$   \\ \hline
$v_1,v_2,v_3,v_6$& $1.05881208405979$   \\ \hline
$v_1,v_2,v_3,v_7$& $0.61626814533203$   \\ \hline

$v_1,v_2,v_4,v_5$& $1.36608045307310$   \\ \hline
$v_1,v_2,v_4,v_6$& $1.06274475848552$   \\ \hline
$v_1,v_2,v_4,v_7$& $0.61717041892812$   \\ \hline

$v_1,v_2,v_5,v_7$& $1.21717866219720$   \\ \hline
$v_1,v_2,v_5,v_7$& $0.64658016206138$   \\ \hline

$v_1,v_2,v_6,v_7$& $0.70102635018388$   \\ \hline
\end{tabular}
\end{table}
\end{nonsec}

\section{Examples: exterior modulus of quadrilaterals}

In this section, we consider several numerical examples to illustrate the accuracy of the proposed method for computing the exterior modulus of quadrilaterals. In the first example, the exact value of the exterior modulus is known. For the second example, we compare the above proposed method against the methods presented in~\cite{HRV2}. 

\begin{nonsec}{\bf Rectangle.}
In this example, we consider the rectangle $R$ with the vertices $z_1=0$, $z_2=\i b$, $z_3=a+\i b$, $z_4=a$ with $a,b>0$. The exact value of the exterior modulus of the quadrilateral $(R; z_1,z_2,z_3,z_4)$ is given by~\cite[p.~82]{ps10}
\begin{equation}\label{eq:emod-R}
{\rm mod}(R; z_1,z_2,z_3,z_4) = \frac{1}{\pi}\,\mu(\kappa),  
\end{equation}
where
\[
\kappa=\psi^{-1}\left(\frac{a}{b}\right), \quad
\psi(\kappa)=\frac{2(\E(\kappa)-(1-\kappa)\K(\kappa))}{\E'(\kappa)-\kappa\K'(\kappa)}.
\]
In our numerical examples below, we assume that $a=1$ and we choose several values of $\kappa$ such that $0.02787\le\kappa\le0.7306$, then $0.02<b=1/\psi(\kappa)<10$. For these values of $b$, the proposed method with $n=2^{13}$ is used to compute approximate values of the exterior modulus of the quadrilateral $(R; z_1,z_2,z_3,z_4)$. The relative error in the computed values is presented in Figure~\ref{fig:emod-Rec} where the exact values of the exterior modulus is computed by~\eqref{eq:emod-R}.

\begin{figure}[H] %
\centerline{
\scalebox{0.7}{\includegraphics[trim=0cm 0cm 0cm 0cm,clip]{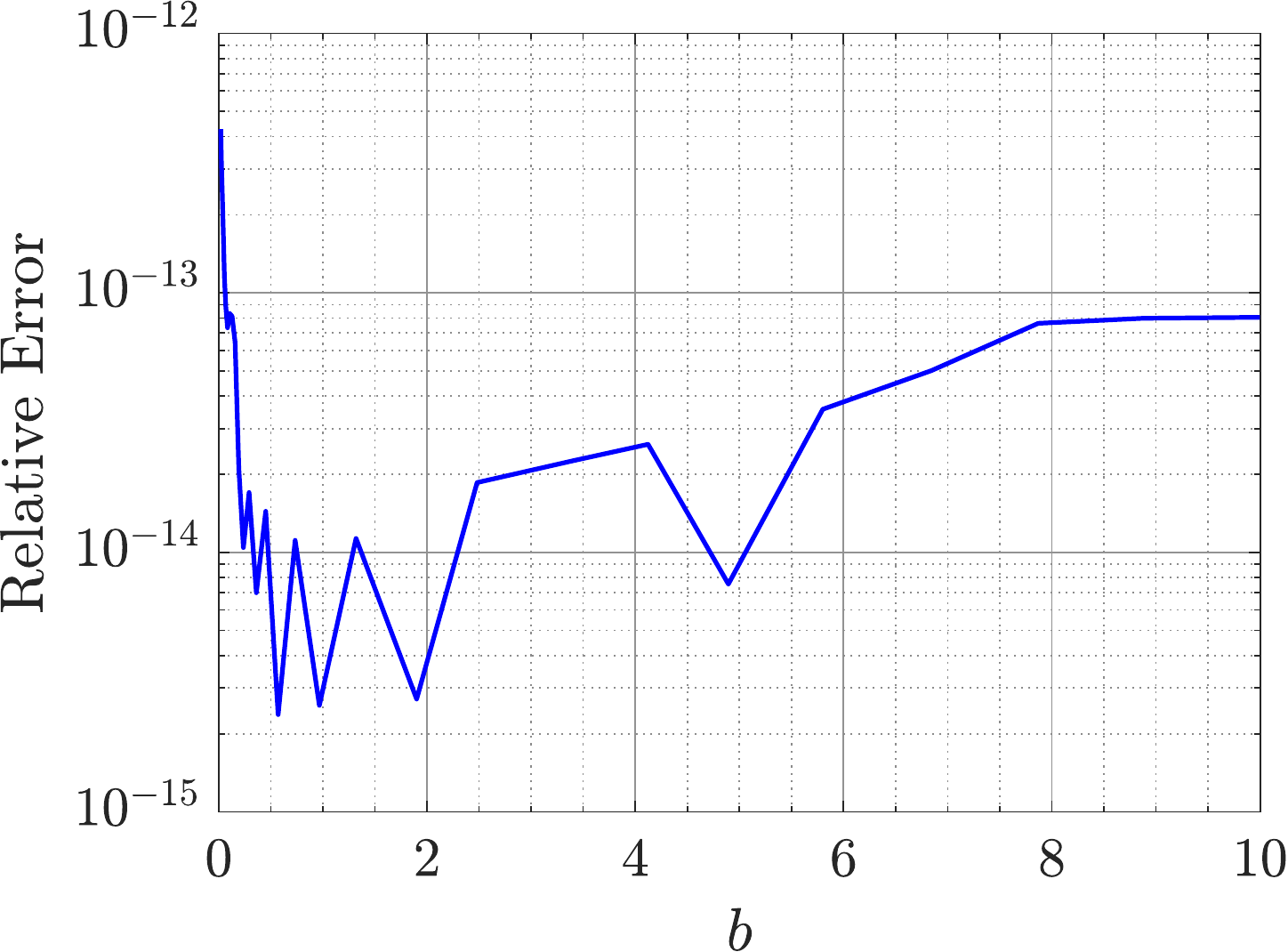}}
}
\caption{The relative error in the computed approximate values of the exterior modulus vs $b$ for the rectangle with the vertices $z_1=0$, $z_2=\i b$, $z_3=1+\i b$, $z_4=1$.}
\label{fig:emod-Rec}
\end{figure}

\end{nonsec}

\begin{nonsec}{\bf Two polygonal quadrilaterals.}
In this example, we compute the exterior modulus for two polygonal quadrilaterals from~\cite{HRV2} (see Figure~\ref{fig:emod-2q}). For the first quadrilateral, we consider the polygon $P_1$ with the vertices $z_1=0$, $z_2=-19/25+\i 21/25$, $z_3=28/25+\i 69/50$, $z_4=1$. In the second quadrilateral, we consider the polygon $P_2$ with the vertices $z_1=0$, $z_2=-3/25+\i 21/25$, $z_3=42/25+4\i$, $z_4=1$. 
The exact value of the exterior modulus of the quadrilaterals $(P_1; z_1,z_2,z_3,z_4)$ and  $(P_2; z_1,z_2,z_3,z_4)$ are unknown. The approximate values of the exterior modulus for these two quadrilaterals are given in Table~\ref{tab:emod-2q} obtained with the proposed method with $n=2^{13}$. Table~\ref{tab:emod-2q} presents also the values of the exterior modulus obtained by three methods presented in~\cite{HRV2} and also the values computed by the SC toolbox in~\cite{HRV2}.

\begin{figure}[H] %
\centerline{
\scalebox{0.55}{\includegraphics[trim=0cm 0cm 0cm 0cm,clip]{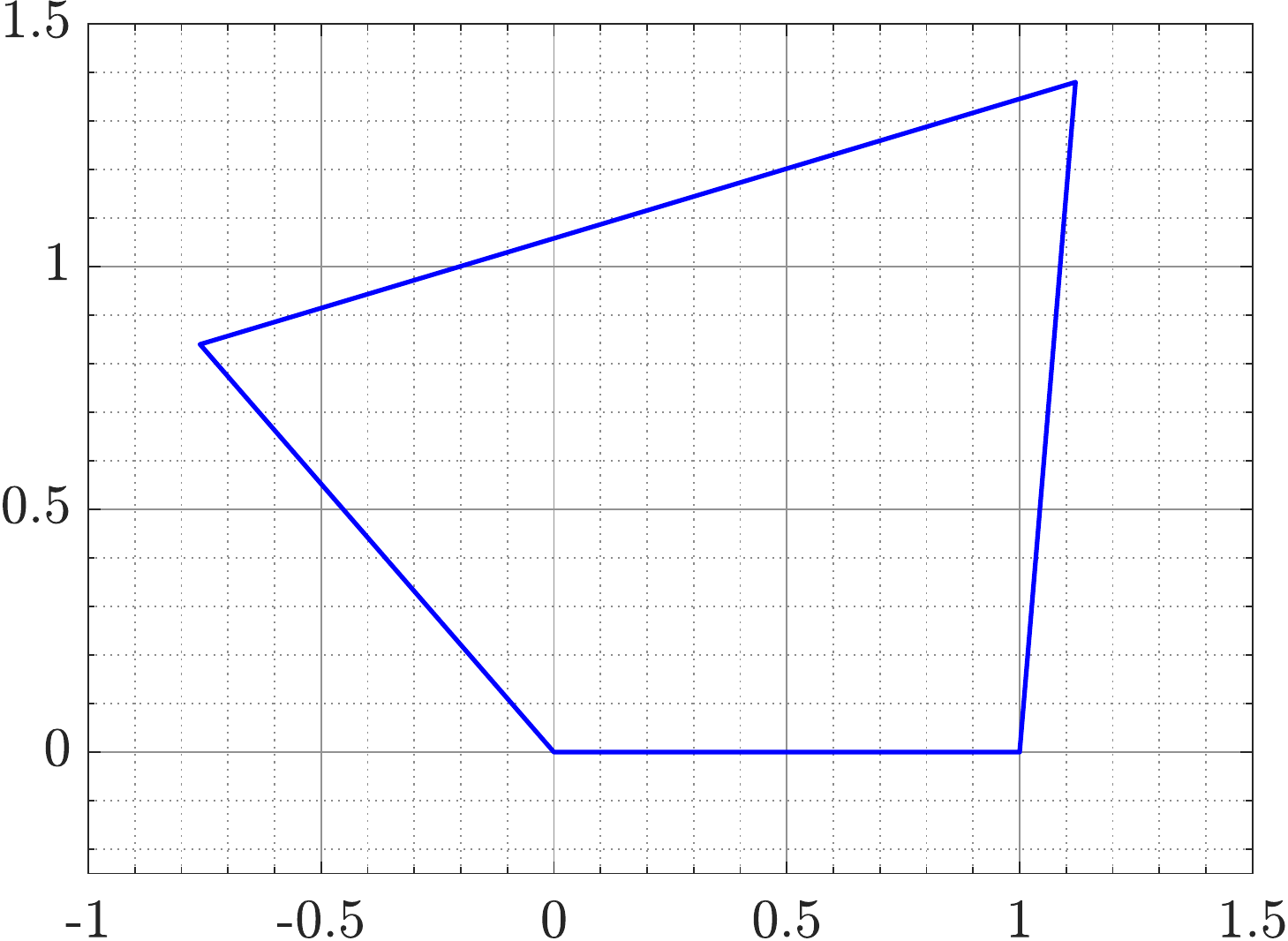}}
\hfill
\scalebox{0.55}{\includegraphics[trim=0cm 0cm 0cm 0cm,clip]{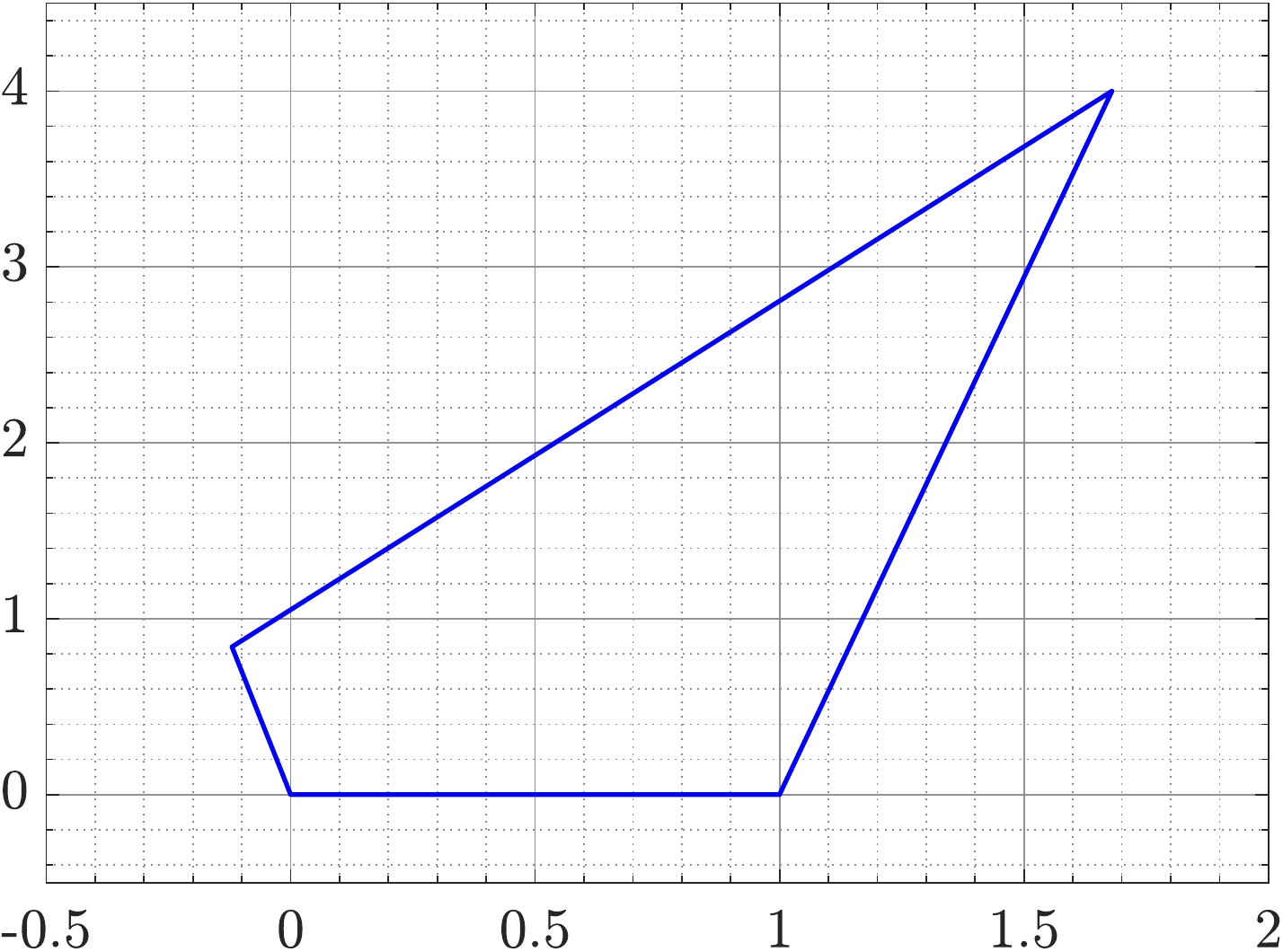}}
}
\caption{The polygon $P_1$ (left) and the polygon $P_2$ (right).}
\label{fig:emod-2q}
\end{figure}

\begin{table}[H]
  \caption{The values of the exterior modulus for the two quadrilaterals $(P_1;z_1, z_2, z_3, z_4)$ and $(P_2;z_1, z_2, z_3, z_4)$.}
\label{tab:emod-2q}  
 	\vspace{0.2cm}
  \begin{tabular}{|l|l|l|} 
 \hline
Method                          & $(P_1;z_1, z_2, z_3, z_4)$ & $(P_2;z_1, z_2, z_3, z_4)$  \\ \hline
Our Method                      & $0.9923416331$  & $0.9592571729$\\ \hline
SC Toolbox~\cite{HRV2}          & $0.9923416323$  & $0.9592571721$ \\ \hline
AFEM~\cite{HRV2}                & $0.9923500126$  & $0.9593012739$ \\ \hline
$hp$-FEM (Interior)~\cite{HRV2} & $0.9923416332$  & $0.9592571731$ \\ \hline
$hp$-FEM (Exterior)~\cite{HRV2} & $0.9923416332$  & $0.9592572007$ \\ \hline
\end{tabular}
\end{table}

\end{nonsec}

\section{Conformal mapping onto gear domains}

\begin{nonsec}{\bf Gear domains.} 
A gear domain $D$ is a special case of the circular arc polygons described above. It is a starlike simply connected domain containing the origin and bordered by arcs of circles centered at the origin and segments of lines passing through the origin~\cite{br,bp,bp2,pe}. Here, we assume that $D$ is bounded. 
The method presented in \S~\ref{sec:map} can be used to compute the conformal mapping $w=\Phi(z)$ for the gear domain $D$ onto the unit disk $\DD$ normalized by~\eqref{eq:Phi-cond} as well as its inverse $z=\Phi^{-1}(w)$ from $\DD$ onto $D$. Assume that $D$ has $m$ vertices $v_k$, $k=1,2,\ldots,m$.
Then the method can be used to compute also the preimages $w_k$, $k=1,2,\ldots,m$, of these vertices.

\end{nonsec}

\begin{nonsec}{\bf A gear domain with $6$ vertices.}
As our first example, we consider a gear domain with $6$ vertices. This example has been considered in~\cite[Fig. 4(b)]{pe} (although the vertices of this domain are not given explicitly in~\cite{pe}, we approximate these vertices from Fig. 4(b) in~\cite{pe}). 
The vertices are given in Table~\ref{tab:g6}.

\begin{table}[H]
  \caption{The vertices $v_k$ and the preimages $w_k$ for the gear domain with $6$ vertices.}
\label{tab:g6}  
 	\vspace{0.2cm}
  \begin{tabular}{|c|c|c|} 
 \hline
$k$   & $v_k$               & $w_k$  \\ \hline
$1$   & $0.75e^{\i\pi/5}$   & $ 0.97953567010215+0.20127064117138\i$       \\ \hline
$2$   & $0.75e^{3\i\pi/5}$  & $ 0.92181215666441+0.38763687624595\i$       \\ \hline
$3$   & $1.25e^{3\i\pi/5}$  & $-0.60224821653432+0.79830889114505\i$       \\ \hline
$4$   & $1.25e^{3\i\pi/2}$  & $-0.82674608972861+0.56257524218406\i$       \\ \hline
$5$   & $e^{3\i\pi/2}$      & $-0.33447524422818-0.94240453680917\i$       \\ \hline
$6$   & $e^{\i\pi/5}$       & $-0.20794689838982-0.97814011647108\i$       \\ \hline
\end{tabular}
\end{table}

The method presented in \S~\ref{sec:map} is used with $n=3\times2^{11}$ to compute the conformal mapping from the gear domain $D$ onto the unit disk and its inverse. 
Figure~\ref{fig:g6} (left) shows the images of several circles $|w|=r$, for $r=0.3,0.45,0.6,0.75,0.9$, in the unit disk under the inverse conformal mapping $z=\Phi^{-1}(w)$.  
The image of the circle $|z|=r$ or part of the circle for $r=0.19,0.58,0.88,1.13$, in the domain $D$ under the conformal mapping $w=\Phi(z)$ is shown in Figure~\ref{fig:g6} (right). The square markers on the unit circle are the preimages of the vertices of the gear domain. The approximate values of the preimages are presented in Table~\ref{tab:g6}.

\begin{figure}[H] %
\centerline{
\scalebox{0.7}{\includegraphics[trim=0cm 0cm 0cm 0cm,clip]{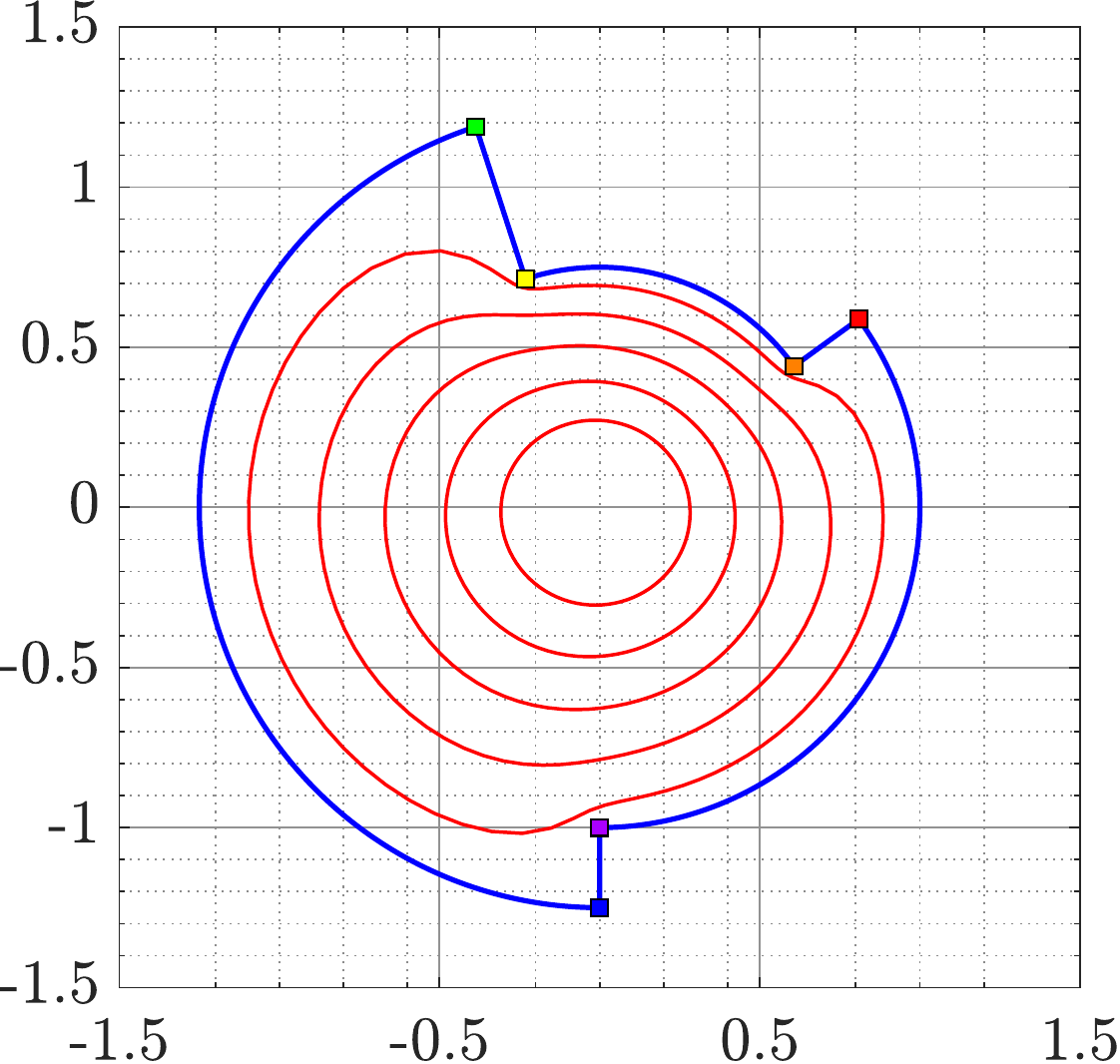}}
\hfill
\scalebox{0.7}{\includegraphics[trim=0cm 0cm 0cm 0cm,clip]{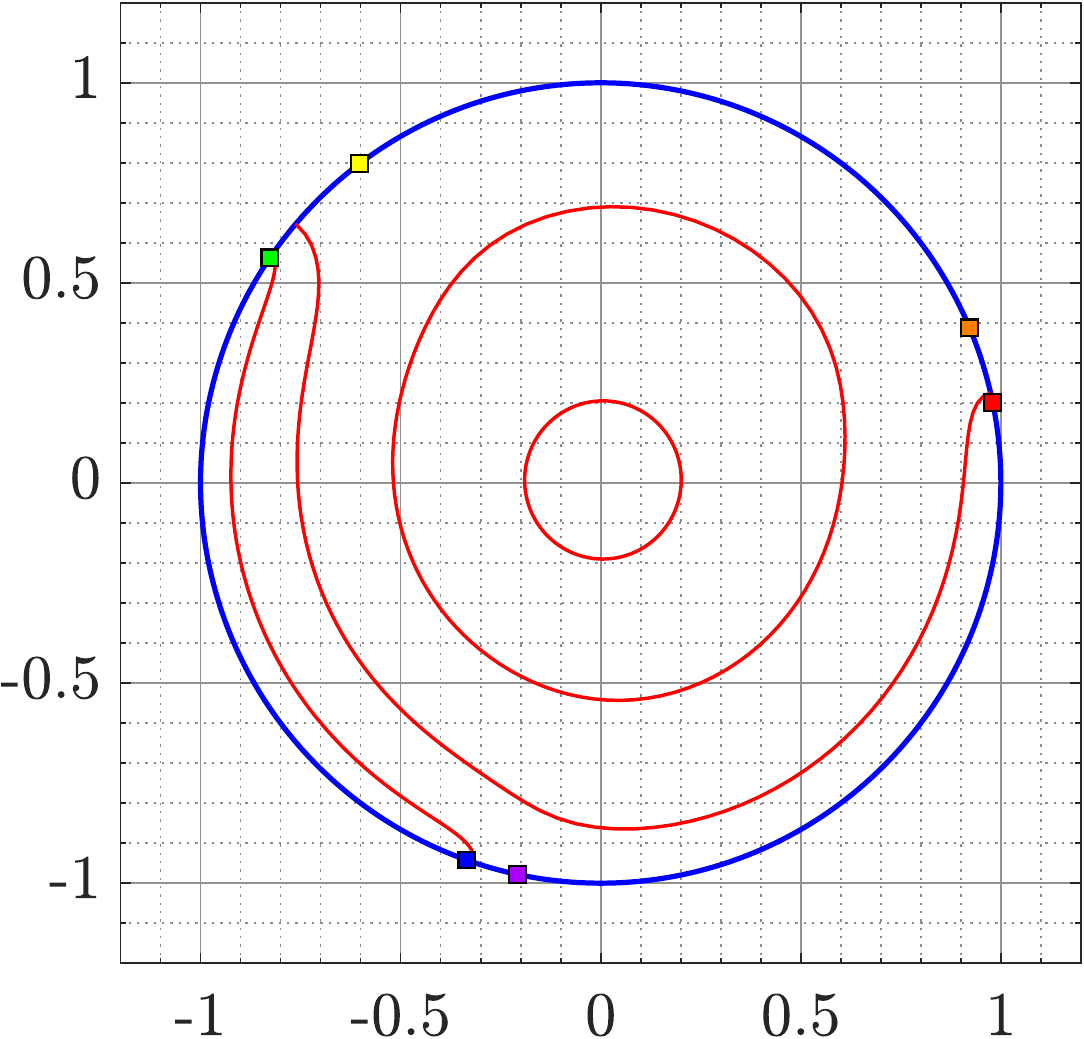}}
}
\caption{The gear domain with $6$ vertices (left) and the unit disk (right).}
\label{fig:g6}
\end{figure}

\end{nonsec}

\begin{nonsec}{\bf One-tooth gear domain.}
Numerical Conformal Mappings onto one-tooth gear domains have been considered in~\cite{br,bp,bp2}. Without loss of generality, a one-tooth gear domain is a circular arc polygonal domain with the vertices $v_1=e^{-\i\theta}$, $v_2=\beta e^{-\i\theta}$, $v_3=\beta e^{\i\theta}$, and $v_4=e^{\i\theta}$ where $0<\theta<\pi$ is the gear angle and $\beta>1$ is the gear ratio~\cite{bp}. 
The method presented in \S~\ref{sec:map} is used with $n=2^{13}$ to compute the conformal mapping from the gear domain $D$ onto the unit disk and its inverse for $\beta=1.5$ and $\theta=\pi/6$. 
Figure~\ref{fig:g1t} (left) shows the images of several circles $|w|=r$, for $r=0.09,0.19,\ldots,0.99$, under the inverse conformal mapping $z=\Phi^{-1}(w)$. The image of the circle $|z|=r$ or part of the circle for $r=0.09,0.19,\ldots,1.49$, under the conformal mapping $w=\Phi(z)$ is shown in Figure~\ref{fig:g6} (right). The square markers on the unit circle indicate the preimages of the vertices of the gear domain. 

\begin{figure}[H] %
\centerline{
\scalebox{0.7}{\includegraphics[trim=0cm 0cm 0cm 0cm,clip]{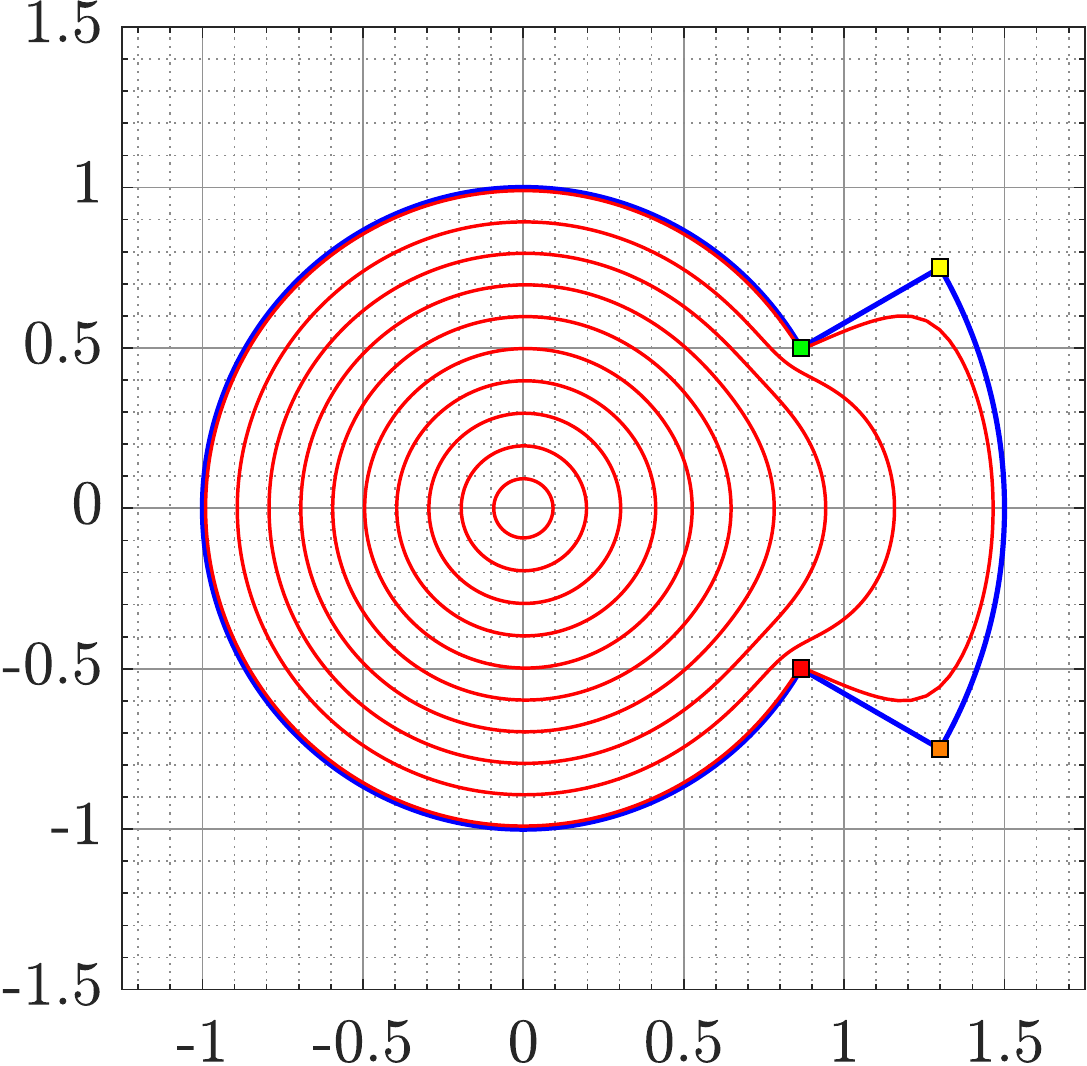}}
\hfill
\scalebox{0.7}{\includegraphics[trim=0cm 0cm 0cm 0cm,clip]{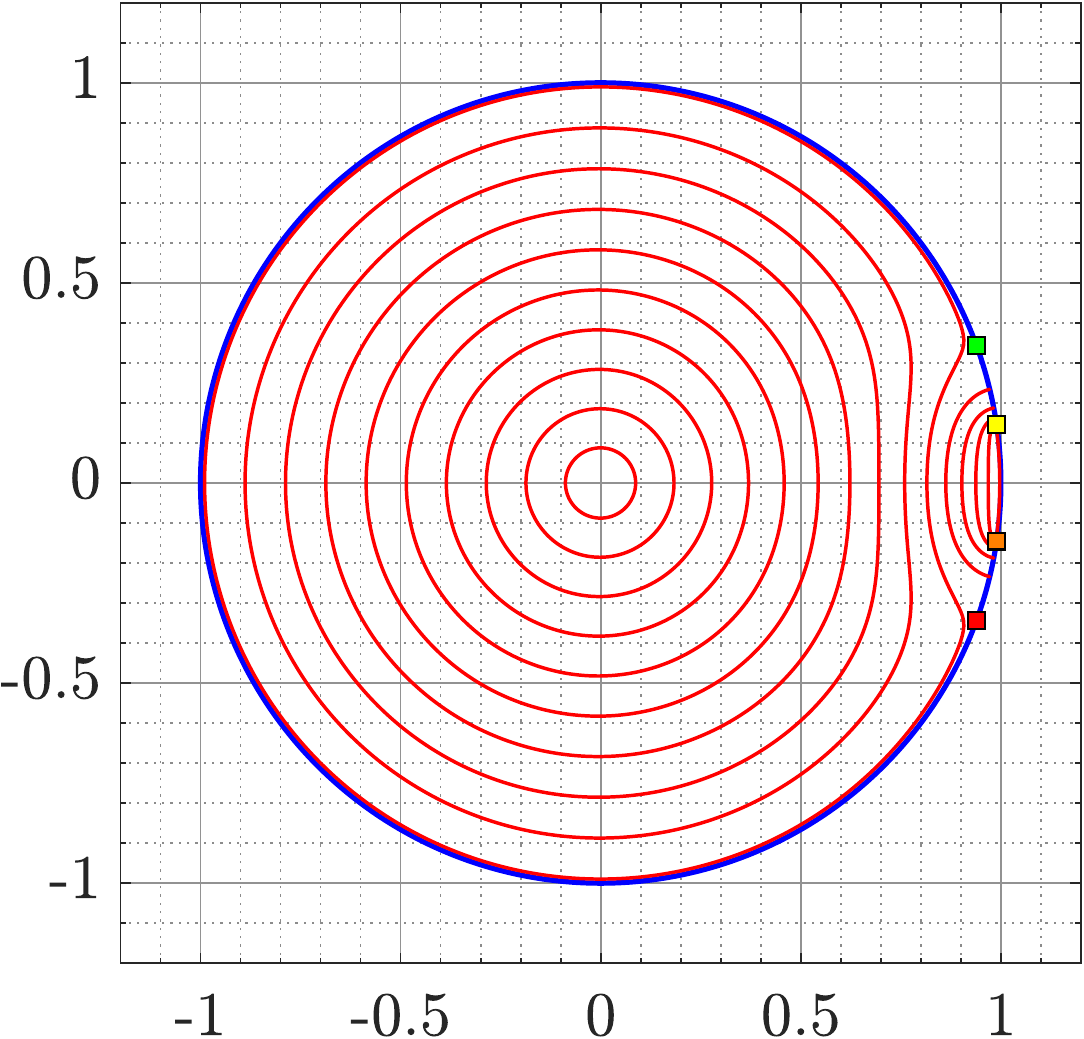}}
}
\caption{The gear domain with one-tooth for $\theta=\pi/6$ and $\beta=1.5$ (left) and the unit disk (right).}
\label{fig:g1t}
\end{figure}

For a fixed $\beta$ and for $0<\theta<\pi$, the modulus of the quadrilateral $(D;e^{-\i\theta},\beta e^{-\i\theta},\beta e^{\i\theta},e^{\i\theta})$ has been computed using the proposed method with $n=2^{13}$. The results for $\beta=1.1,1.25,1.5,2$ are presented in Figure~\ref{fig:g1tm}. It is clear from Figure~\ref{fig:g1tm} that the modulus approaches zero as $\theta\to0$ or $\theta\to\pi$. Further, the results presented in Figure~\ref{fig:g1tm} validate numerically the conjuncture in~\cite[p.~90]{bp2}. That is, there are exactly two gears corresponding to two different values of $\theta$ with the same modulus except for one value of $\theta\in(0,\pi)$ where the modulus has its maximum value. These maximums are marked with squares in Figure~\ref{fig:g1tm}. The location of the maximum value depends on the value of $\beta$ and it moves towards $\pi$ as $\beta$ increases. Moreover, the maximum value of the modulus increases as $\beta$ decreases toward $1$.

\begin{figure}[H] %
\centerline{
\scalebox{0.7}{\includegraphics[trim=0cm 0cm 0cm 0cm,clip]{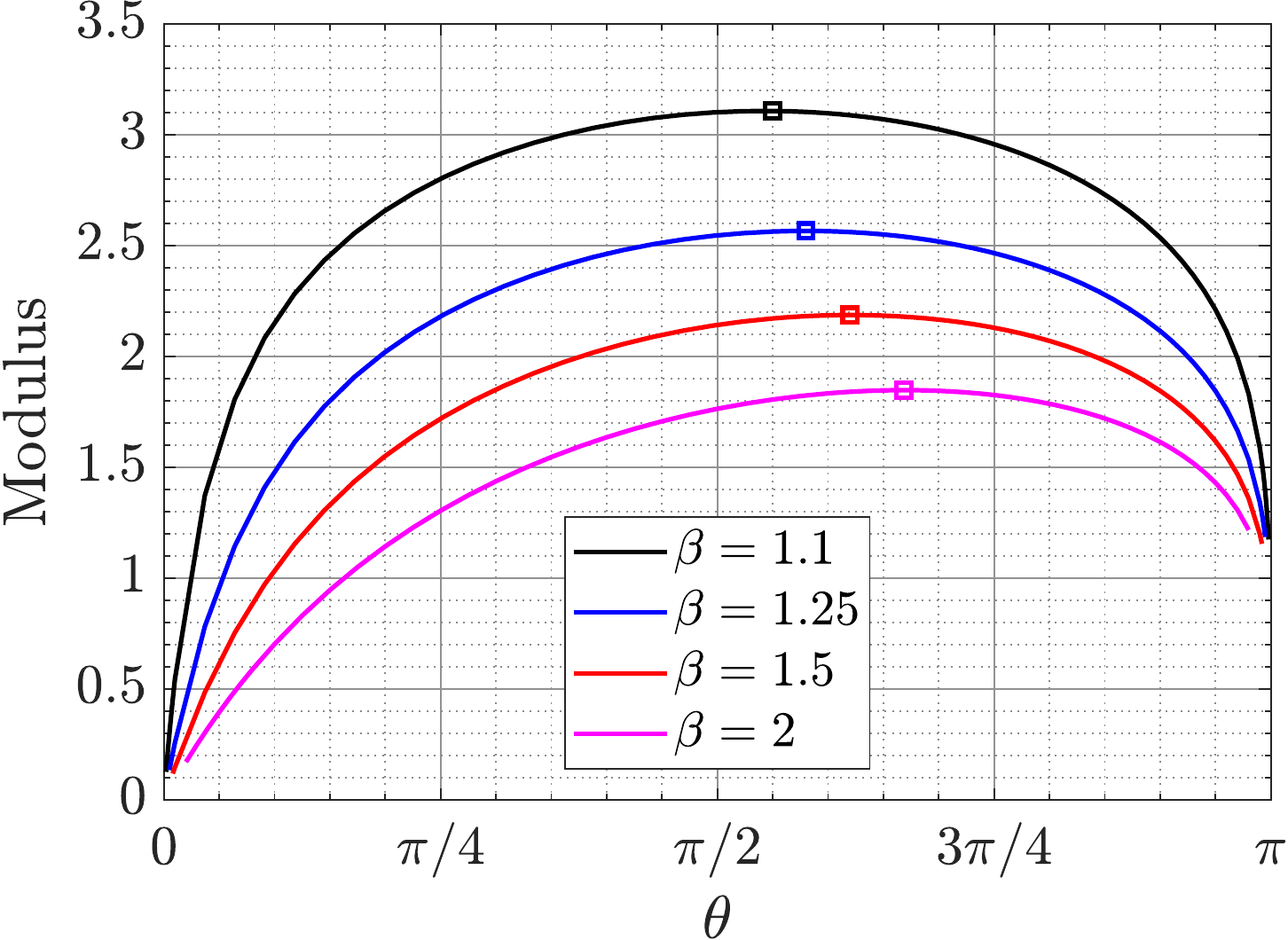}}
}
\caption{The modulus of the quadrilateral $(D;e^{-\i\theta},\beta e^{-\i\theta},\beta e^{\i\theta},e^{\i\theta})$.}
\label{fig:g1tm}
\end{figure}

\end{nonsec}

\begin{nonsec}{\bf A multitooth gear domain.}
We consider a multitooth gear domain with $12$ vertices as in Table~\ref{tab:g12t}.
The method presented in \S~\ref{sec:map} is used with $n=3\times2^{12}$ to compute the conformal mapping from the gear domain $D$ onto the unit disk and its inverse.
Figure~\ref{fig:g12t} (left) shows the images of several circles $|w|=r$, for $r=0.09,0.19,\ldots,0.99$, under the inverse conformal mapping $z=\Phi^{-1}(w)$. The image of the circle $|z|=r$ or part of the circle for $r=0.09,0.19,\ldots,1.99$, under the conformal mapping $w=\Phi(z)$ is shown in Figure~\ref{fig:g12t} (right). 

\begin{table}[H]
  \caption{The vertices $v_k$ and the preimages $w_k$ for the multitooth gear domain with $12$ vertices.}
\label{tab:g12t}  
 	\vspace{0.2cm}
  \begin{tabular}{|c|c|c|} 
 \hline
$k$   & $v_k$               & $w_k$  \\ \hline
$1$   & $e^{\i\pi/6}$       & $ 0.86701428817497+0.49828327696246\i$       \\ \hline
$2$   & $e^{\i\pi/2}$       & $-0.28316473230969+0.95907128743174\i$       \\ \hline
$3$   & $2e^{\i\pi/2}$      & $-0.56900711726358+0.82233259725210\i$       \\ \hline
$4$   & $2e^{3\i\pi/4}$     & $-0.65069062054555+0.75934295040781\i$       \\ \hline
$5$   & $1.5e^{3\i\pi/4}$   & $-0.71186505065025+0.70231627466742\i$       \\ \hline
$6$   & $1.5e^{\i\pi}$      & $-0.95549393111898+0.29501075843909\i$       \\ \hline
$7$   & $1.25e^{\i\pi}$     & $-0.97908358634907+0.20345842558580\i$       \\ \hline
$8$   & $1.25e^{3\i\pi/2}$  & $-0.62508676492722-0.78055527434822\i$       \\ \hline
$9$   & $0.75e^{3\i\pi/2}$  & $-0.32775376595675-0.94476318138524\i$       \\ \hline
$10$  & $0.75e^{11\i\pi/6}$ & $ 0.97086850902193-0.23961289236922\i$       \\ \hline
$11$  & $1.75e^{11\i\pi/6}$ & $ 0.98506920087238+0.17215884959145\i$       \\ \hline
$12$  & $1.75e^{\i\pi/6}$   & $ 0.95294901093885+0.30313063611365\i$       \\ \hline
\end{tabular}
\end{table}

\begin{figure}[H] %
\centerline{
\scalebox{0.7}{\includegraphics[trim=0cm 0cm 0cm 0cm,clip]{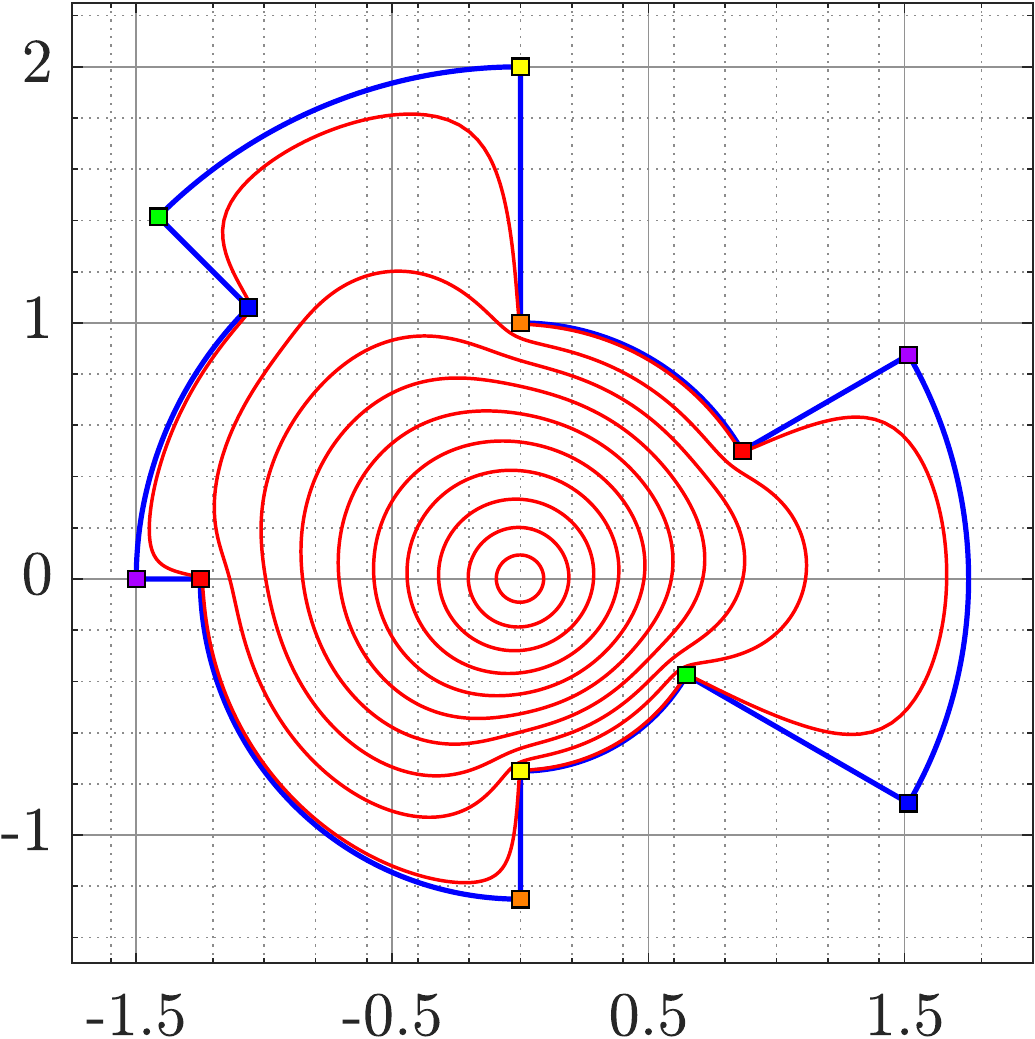}}
\hfill
\scalebox{0.7}{\includegraphics[trim=0cm 0cm 0cm 0cm,clip]{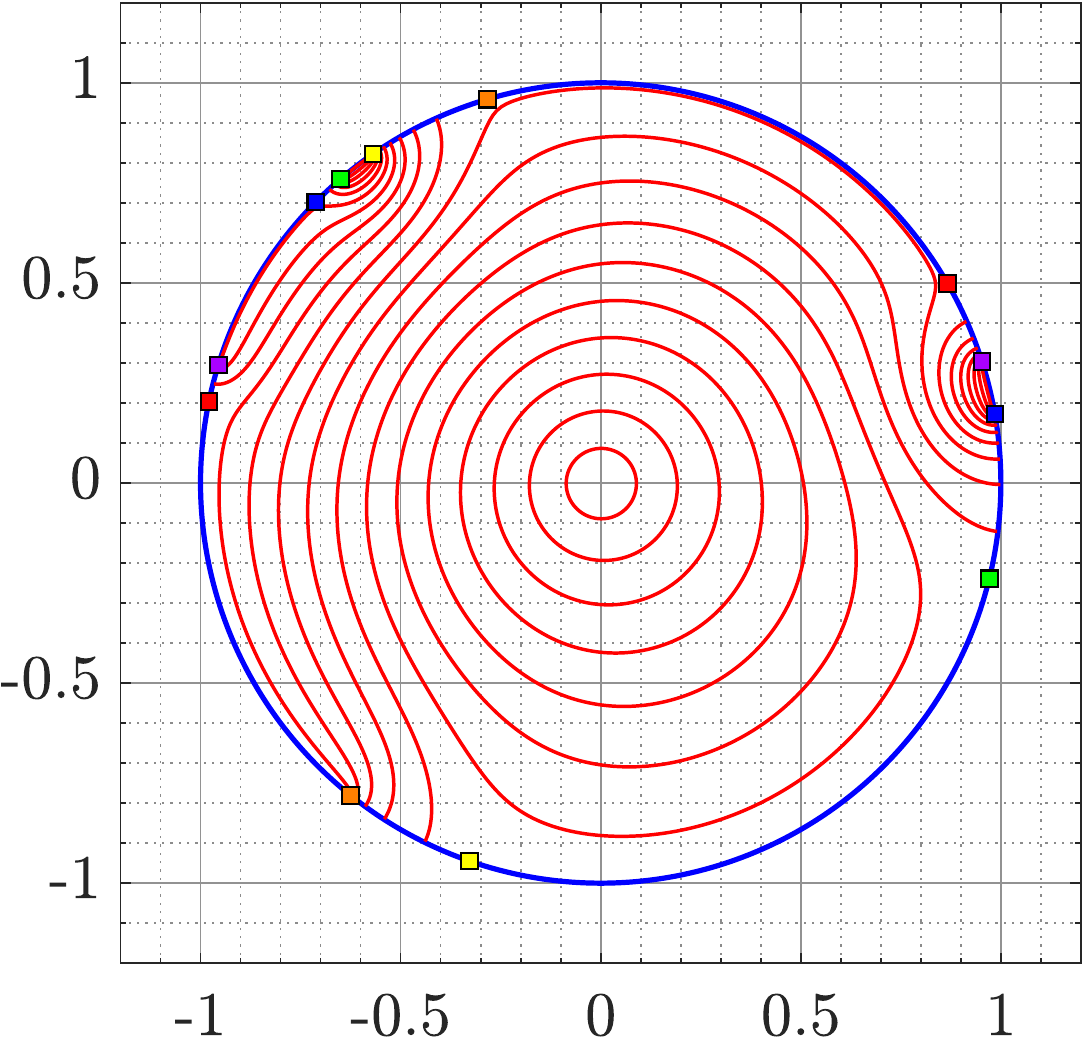}}
}
\caption{The multitooth gear domain (left) and the unit disk (right).}
\label{fig:g12t}
\end{figure}

\end{nonsec}

\begin{nonsec}{\bf Complement of an Annular Rectangle.}
Consider the circular arc polygon whose boundary consists of the straight segment from $e^{-\i\theta}$ to $\beta^2 e^{-\i\theta}$, the circular arc from $\beta^2 e^{-\i\theta}$ to  $\beta^2 e^{\i\theta}$, the straight segment from $\beta^2 e^{\i\theta}$ to $e^{\i\theta}$, and the circular arc from $e^{\i\theta}$ to $e^{-\i\theta}$ where $\theta\in(0,\pi)$ and $\beta>1$ (see domain $G$ in Figure~\ref{fig:compl} for $\theta=\pi/4$ and $\beta=1.5$). This domain $G$ is called an \emph{annular rectangle}~\cite{bp}. Consider also the gear domain $D$ with the vertices $e^{-\i\theta}$, $\beta e^{-\i\theta}$, $\beta e^{\i\theta}$, and $e^{\i\theta}$ (see domain $D$ in Figure~\ref{fig:compl}). Then, it follows from~\cite[Theorem 6.1]{bp} that the exterior modulus of the quadrilateral $(G;e^{-\i\theta},\beta^2 e^{-\i\theta},\beta^2 e^{\i\theta},e^{\i\theta})$ is half of ${\rm mod}(D;e^{-\i\theta},\beta e^{-\i\theta},\beta e^{\i\theta},e^{\i\theta})$.

\begin{figure}[H] %
\centerline{
\scalebox{0.7}{\includegraphics[trim=0cm 0cm 0cm 0cm,clip]{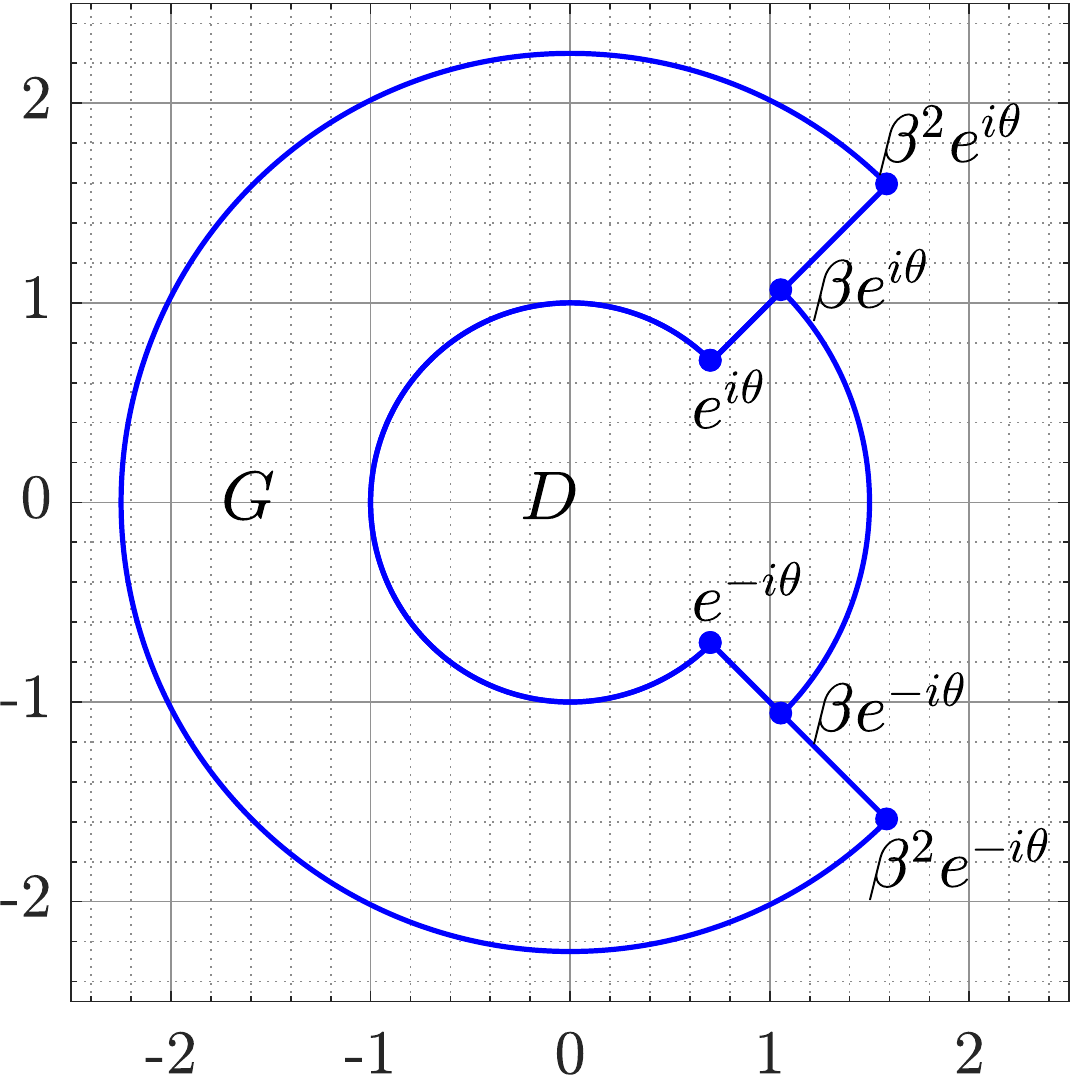}}
}
\caption{The annular rectangle domain $G$ and the gear domain $D$.}
\label{fig:compl}
\end{figure}

In this example, we use our proposed method with $n=2^{13}$ to compute the exterior modulus of the quadrilateral $(G;e^{-\i\theta},\beta^2 e^{-\i\theta},\beta^2 e^{\i\theta},e^{\i\theta})$ as well as ${\rm mod}(D;e^{-\i\theta},\beta e^{-\i\theta},\beta e^{\i\theta},e^{\i\theta})$ for several values of $\theta$. The absolute value of the difference between the computed exterior modulus and $0.5\,{\rm mod}(D;e^{-\i\theta},\beta e^{-\i\theta},\beta e^{\i\theta},e^{\i\theta})$ is considered as the error in the computed values. The obtained results are presented in Table~\ref{tab:compl}.

\begin{table}[H]
  \caption{The exterior modulus of the quadrilateral $(G;e^{-\i\theta},\beta^2 e^{-\i\theta},\beta^2 e^{\i\theta},e^{\i\theta})$ and ${\rm mod}(D;e^{-\i\theta},\beta e^{-\i\theta},\beta e^{\i\theta},e^{\i\theta})$.}
\label{tab:compl}  
 	\vspace{0.2cm}
  \begin{tabular}{|c|c|c|c|} 
 \hline
$\theta$   & Exterior modulus      & $0.5\,{\rm mod}(D;e^{-\i\theta},\beta e^{-\i\theta},\beta e^{\i\theta},e^{\i\theta})$  & Error  \\ \hline
$0.1\pi$ & $0.51830606688359$  & $0.51830606688379$  & $1.96\times10^{-13}$ \\ \hline
$0.2\pi$ & $0.77581840983574$  & $0.77581840983561$  & $1.30\times10^{-13}$ \\ \hline
$0.3\pi$ & $0.92576131108263$  & $0.92576131108211$  & $5.21\times10^{-13}$ \\ \hline
$0.4\pi$ & $1.01795618251692$  & $1.01795618251687$  & $4.91\times10^{-14}$ \\ \hline
$0.5\pi$ & $1.07133752300218$  & $1.07133752300216$  & $1.40\times10^{-14}$ \\ \hline
$0.6\pi$ & $1.09298754180547$  & $1.09298754180560$  & $1.27\times10^{-13}$ \\ \hline
$0.7\pi$ & $1.08332598419075$  & $1.08332598419105$  & $2.95\times10^{-13}$ \\ \hline
$0.8\pi$ & $1.03535729272694$  & $1.03535729272677$  & $1.68\times10^{-13}$ \\ \hline
$0.9\pi$ & $0.92271712140416$  & $0.92271712140442$  & $2.59\times10^{-13}$ \\ \hline
\end{tabular}
\end{table}

\end{nonsec}

\newpage
\def\cprime{$'$} \def\cprime{$'$} \def\cprime{$'$}
\providecommand{\bysame}{\leavevmode\hbox to3em{\hrulefill}\thinspace}
\providecommand{\MR}{\relax\ifhmode\unskip\space\fi MR }
\providecommand{\MRhref}[2]{%
  \href{http://www.ams.org/mathscinet-getitem?mr=#1}{#2}
}
\providecommand{\href}[2]{#2}

\end{document}